\documentclass[12pt,reqno]{amsart}

\usepackage[left=2cm,right=2cm,top=2cm,bottom=2cm]{geometry}
\usepackage{amsmath, amsthm, amscd, amsfonts, amssymb, graphicx, color, mathrsfs}
\usepackage[bookmarksnumbered, colorlinks, plainpages]{hyperref}
\usepackage{cite}
\usepackage[all]{xy}
\usepackage{slashed}
\usepackage{tikz-cd}
\usepackage{mathabx}
\usepackage{tipa}
\usepackage{soul}
\usepackage{cancel}
\usepackage{ulem}
\usepackage{algorithm}
\usepackage{algpseudocode}
\usepackage{enumitem}
\usepackage{multirow}
\usepackage{float}
\usepackage{verbatim}
\usepackage{graphicx}

\allowdisplaybreaks

\newlength\tindent
\newtheorem{theorem}{Theorem}[section]
\newtheorem{lemma}[theorem]{Lemma}

\newtheorem{proposition}[theorem]{Proposition}
\newtheorem{corollary}[theorem]{Corollary}

\theoremstyle{definition}

\newtheorem{example}[theorem]{Example}

\theoremstyle{remark}
\newtheorem{remark}[theorem]{Remark}
\numberwithin{equation}{section}

\begin{document}
\setcounter{page}{1}

\title[Balanced non-Kähler metrics on complex projective spaces]{Inhomogeneous Balanced non-K\"ahler metrics on complex projective spaces}

\author[B. Grajales]{Brian Grajales}
\address{Brian Grajales \endgraf
Universidade Estadual de Maringá, Departamento de Matem\'{a}tica, Avenida Colombo, 5790, Campus Universitário, 87020-900, Maringá - PR, Brazil 
\endgraf
  {\it E-mail address:} {\rm bdgtriana@uem.br}
  }

\author[L. Grama]{Lino Grama}
\address{Lino Grama \endgraf
IMECC-Unicamp, Departamento de Matem\'{a}tica. Rua S\'{e}rgio Buarque de Holanda \endgraf 651, Cidade Universit\'{a}ria Zeferino Vaz. 13083-859, Campinas - SP, Brazil
\endgraf
{\it E-mail address:} {\rm lino@ime.unicamp.br}
  }



\begin{abstract}
We consider the cohomogeneity one action of $G_{p,q}=\mathrm{U}(p+1)\times\mathrm{U}(q+1)$ on $\mathbb{C}P^{p+q+1}$, where $p\geq q\geq 1$ and $p>1$. We characterize all $G_{p,q}$-invariant Riemannian metrics compatible with the canonical complex structure and obtain parametrizations of the invariant K\"ahler and balanced metrics. As an application, we use these parametrizations to construct a family of inhomogeneous balanced non-K\"ahler metrics with positive scalar curvature on complex projective spaces of odd complex dimension.
\end{abstract}

\maketitle

\allowdisplaybreaks

\section{Introduction}
Let $(M,J)$ be a complex manifold of complex dimension $n\geq 2$, and let $g$ be a Hermitian metric on $(M,J)$ with fundamental form $\omega.$ The metric $g$ is called K\"ahler if $d\omega=0$, and balanced if $d\omega^{n-1}=0$. The terminology {\it balanced} was introduced by Michelsohn \cite{M1982}, who studied these metrics and characterized compact complex manifolds admitting them. Every K\"ahler metric is balanced, but the converse is false for $n\geq 3$. Examples of balanced non-K\"ahler metrics have been constructed on several classes of compact complex manifolds. Fu, Li, and Yau \cite{FLY2012} constructed such metrics on non-K\"ahler Calabi-Yau threefolds. Correa \cite{C2023} studied special Hermitian metrics on compact non-K\"ahler Calabi-Yau manifolds and obtained further examples of balanced metrics. For a broader discussion of balanced metrics in complex geometry, we refer to Fino \cite{F2024}.\\

Lie group symmetries provide a natural framework for studying balanced Hermitian metrics. Two relevant settings are complex homogeneous spaces, on which a Lie group acts transitively by biholomorphisms, and cohomogeneity one manifolds, whose principal orbits have codimension one. In the homogeneous setting, Fino, Grantcharov, and Vezzoni \cite{FGV2019} obtained necessary and sufficient conditions for the existence of balanced metrics on compact complex homogeneous spaces admitting an invariant volume form. In particular, they showed that every invariant Hermitian metric on a generalized flag manifold is balanced. More recently, Lauret \cite{L2026} investigated invariant balanced Hermitian metrics on compact homogeneous complex manifolds and established further existence results. Related questions for left-invariant Hermitian structures on Lie groups and Lie algebras have been studied by Andrada and Villacampa \cite{AV2016} and Fino and Paradiso \cite{FP2023}. In contrast to this flexibility in the homogeneous setting, the balanced condition may exhibit significant rigidity on cohomogeneity one manifolds. Indeed, Angella and Pediconi \cite{AP2022} proved that the balanced and K\"ahler conditions are equivalent within their class of standard B\'erard--Bergery cohomogeneity one Hermitian manifolds. Alonso and Salvatore \cite{AS2021,AS2024} investigated invariant balanced non-K\"ahler $\mathrm{SU}(3)$-structures on simply connected six-dimensional cohomogeneity one manifolds, including questions concerning their extension across singular orbits.\\

These contrasting phenomena motivate the study of the balanced condition for other cohomogeneity one actions on compact complex manifolds. In particular, it is natural to ask whether such actions admit invariant balanced non-K\"ahler metrics and, if so, whether these metrics can be characterized globally. Complex projective spaces provide a natural setting for this investigation. Their cohomogeneity one actions were classified, up to orbit equivalence, by Takagi \cite{T1973}. These actions have also been considered in the study of invariant Hermitian and Einstein metrics. Angella and Pediconi \cite{AP2022} constructed invariant non-K\"ahler Hermitian metrics conformal to the Fubini--Study metric on $\mathbb{C}P^n$. In a previous work with Araujo \cite{AGG2026}, we studied diagonal invariant Einstein metrics on complex projective spaces under the assumption that the singular orbits are totally geodesic. Broder \cite{B2026} proved that every $\mathrm{SO}(3)$-invariant Einstein metric on $\mathbb{C}P^2$ is homothetic to the Fubini--Study metric. More recently, Cao-Labora and Rodr\'iguez-V\'azquez \cite{CR2026} constructed inhomogeneous Einstein metrics on $\mathbb{C}P^n$ for $n\in\{3,4,5,6,7\}$. Our focus here is instead on the characterization and construction of balanced Hermitian metrics compatible with the canonical complex structure.\\

More precisely, let $p\geq q\geq 1$, with $p>1$, and consider the action of $G_{p,q}:=\mathrm{U}(p+1)\times\mathrm{U}(q+1)$ on $\mathbb{C}P^{p+q+1}$ defined by $(U,V)\cdot[z:w]=[Uz:Vw]$, where $z\in\mathbb{C}^{p+1}$ and $w\in\mathbb{C}^{q+1}$. This action has cohomogeneity one, with two singular orbits diffeomorphic to $\mathbb{C}P^p$ and $\mathbb{C}P^q$. Our main objective is to characterize all $G_{p,q}$-invariant balanced Hermitian metrics on $\mathbb{C}P^{p+q+1}$. We first describe the invariant Riemannian metrics on the principal part and determine which of them are compatible with the canonical complex structure $J$. Using the smoothness criteria of Verdiani and Ziller \cite{VZ}, we obtain necessary and sufficient conditions for these metrics to extend smoothly across both singular orbits. The resulting invariant Hermitian metrics are parametrized by triples $(h_0,h_2,h_3)$ of smooth functions on $[0,1]$ satisfying explicit positivity and boundary conditions (see Proposition~\ref{prop:smoothness-hj}). In terms of this parametrization, the K\"ahler condition reduces to a system of two linear first-order differential equations, whereas the balanced condition reduces to a single first-order differential equation. Solving these equations, we obtain a global parametrization of all $G_{p,q}$-invariant balanced Hermitian metrics by pairs $(h_0,V)$ of smooth functions on $[0,1]$ satisfying $h_0>0$ and $V(0)=0$ (see Corollary~\ref{coro:parametrization-balanced-metrics}). Moreover, the K\"ahler metrics correspond precisely to $V=0$, and hence are parametrized by positive smooth functions $h_0$. Thus, the balanced condition admits a functional family of non-K\"ahler deformations within this invariant setting. As an application, we consider the case $p=q=m$, with $m\geq 2$, and construct an explicit one-parameter family $\{g_\varepsilon\}_{\varepsilon\in\mathbb{R}}$ of $G_{m,m}$-invariant balanced Hermitian metrics on $\mathbb{C}P^{2m+1}$. For this family, $g_0$ is the Fubini--Study metric and $g_\varepsilon$ is non-K\"ahler for every $\varepsilon\neq 0$. Furthermore, we prove that, for sufficiently small nonzero $\varepsilon$, the scalar curvature of $g_\varepsilon$ is positive and non-constant. Consequently, these metrics are inhomogeneous.\\

The paper is organized as follows. In Section~\ref{section:preliminaries}, we recall results on invariant metrics on cohomogeneity one manifolds. In Section~\ref{section:Gpq-action}, we describe the $G_{p,q}$-action on complex projective spaces and characterize the $G_{p,q}$-invariant Hermitian metrics. In Section~\ref{section:balanced-metrics}, we obtain the parametrization of all $G_{p,q}$-invariant balanced Hermitian metrics. Finally, in Section~\ref{section:application}, we construct a family of inhomogeneous balanced non-K\"ahler metrics with positive scalar curvature on complex projective spaces of odd complex dimension.

\section{Preliminaries}\label{section:preliminaries}
\subsection{Invariant metrics in cohomogeneity one}\label{subsec:invariant-metrics}
Let $G$ be a compact connected Lie group, $M$ a connected manifold, and 
\[
\phi:G\times M\longrightarrow M,\ \phi(a,x)=a\cdot x
\]
a smooth cohomogeneity one action. For each $a\in G$, we define
\[
\phi_a:M\longrightarrow M,\ \phi_a(x):=a\cdot x,
\]
which is a diffeomorphism for each $a\in G.$ We assume that the action has no exceptional orbits and that the orbit space $M/G$ is diffeomorphic to the closed interval $[0,b],$ for some $b>0.$ Let
\[
\pi:M\longrightarrow M/G\approx[0,b]
\]
be the orbit projection. The boundary points $0$ and $b$ correspond to the singular orbits and every point in $(0,b)$ corresponds to a principal orbit. Fix a principal isotropy group $H$ and choose a smooth curve $\gamma:[0,b]\longrightarrow M$ such that $\pi(\gamma(t))=t,\ 0\leq t\leq b$ and $G_{\gamma(t)}=H, \ 0<t<b.$ Thus, $\gamma$ intersects every $G$-orbit at exactly one point. Set $Q_0:=G\cdot\gamma(0),\ Q_1:=G\cdot\gamma(b),$ and denote the isotropy group of $\gamma(0),\gamma(b)$ by $K_0:=G_{\gamma(0)},\ K_1:=G_{\gamma(b)},$ respectively. For $t\in(0,b)$, the orbit $P_t:=G\cdot\gamma(t)$ is principal and satisfies $P_t\approx G/H.$ In this setting, we have that $H\subseteq K_j\subseteq G,$ $Q_j\approx G/K_j$ and $K_j/H\approx\mathbb S^{\ell_j},\ j=0,1,$ for some $\ell_j\geq1$. Let $M_0:=M\setminus(Q_0\cup Q_1)$ be the union of all principal orbits. Then, the map
\[
\Psi:(0,b)\times G/H\ni(t,aH)\longmapsto a\cdot \gamma(t)\in M_0,
\]
is a $G$-equivariant diffeomorphism, where we define the action of $G$ on $(0,b)$ to be trivial. In particular, every point $x\in M_0$ can be written as $x=a\cdot \gamma(t),$ for some $t\in(0,b)$ and $a\in G,$ and the parameter $t$ induces a smooth $G$-invariant function $t:M_0\longrightarrow(0,b),$ defined by
\[
t(x)=t_0
\Longleftrightarrow x\in G\cdot\gamma(t_0).
\]
Consequently, $dt$ is a $G$-invariant $1$-form on $M_0$ and it satisfies $dt_{\gamma(t)}(\gamma'(t))=1,$ $dt_{\gamma(t)}(v)=0,$ for every $v\in T_{\gamma(t)}P_t.$\\

Let $\mathfrak{g}, \mathfrak{h}$ and $\mathfrak{k}_j$ denote the Lie algebras of $G,$ $H$ and $K_j$ respectively. For a given $X\in\mathfrak{g},$ we consider the fundamental field $X^*\in\mathfrak{X}(M)$ given by
\[
X^*(x):=\left.\frac{d}{ds}\exp(sX)\cdot x\,\right|_{s=0},
\]
where $\exp:\mathfrak{g}\longrightarrow G$ is the exponential map of $G.$ Let $\partial_t$ denote the vector field on $(0,b)\times G/H$ tangent to the first factor. We define
\[
\widetilde{\gamma'}:=\Psi_*\left(\partial_t\right)\in\mathfrak{X}(M_0).
\]
Since $\Psi(t,eH)=\gamma(t)$, we have $\widetilde{\gamma'}(\gamma(t))=\gamma'(t),\ 0<t<b.$ Fix an $\mathrm{Ad}(G)$-invariant inner product $(\cdot,\cdot)$ on $\mathfrak{g}$ and let $\mathfrak{m}$ be the $(\cdot,\cdot)$-orthogonal complement of $\mathfrak{h}$ in $\mathfrak{g}.$ For each $t\in(0,b),$ we identify $\mathfrak{m}$ with $T_{\gamma(t)}P_t$ via $\mathfrak{m}\ni X\longmapsto X^*(\gamma(t))\in T_{\gamma(t)}P_t.$ Since $\gamma$ is transverse to the $G$-orbits, we have
\[
T_{\gamma(t)}M=\mathbb{R}\gamma'(t)\oplus T_{\gamma(t)}P_t\simeq\mathbb{R}\gamma'(t)\oplus\mathfrak{m},\ 0<t<b.
\]
We consider the family of representations 
\[
\rho_t:H\longrightarrow\mathrm{GL}(T_{\gamma(t)}M),\ \rho_t(h):=(d\phi_h)_{\gamma(t)},\ 0<t<b.
\]
Then, $\rho_t(h)\gamma'(t)=\gamma'(t)$ and $\rho_t(h)(X^*(\gamma(t)))=(\mathrm{Ad}(h)X)^*(\gamma(t)),$ for every $h\in H$ and $X\in\mathfrak{m}.$ 
\begin{proposition}\label{prop:invariant-metrics}
If $g$ is a $G$-invariant Riemannian metric on $M_0$, then, for each $t\in(0,b)$, $g_{\gamma(t)}$ is a $\rho_t(H)$-invariant inner product on $T_{\gamma(t)}M$. Moreover, the family $(g_{\gamma(t)})_{t\in(0,b)}$ is smooth, in the sense that for any basis $\{X_1,...,X_n\}$ of $\mathfrak{m},$ all the functions
\begin{equation}\label{eq:metric-Ej-Ek}
t\longmapsto g_{\gamma(t)}(E_j(t),E_k(t)),\ j,k=0,...,n,
\end{equation}
are smooth, where $E_0(t):=\gamma'(t),$ $E_j(t):=X^*_j(\gamma(t)),$ $j=1,...,n.$\\

Conversely, let $(\langle\cdot,\cdot\rangle_t)_{t\in(0,b)}$ be a family such that, for each $t\in(0,b),$ $\langle\cdot,\cdot\rangle_t$ is a $\rho_t(H)$-invariant inner product on $T_{\gamma(t)}M.$ Assume that this family is smooth in the sense above. Then there exists a unique smooth $G$-invariant Riemannian metric $g$ on $M_0$ such that $g_{\gamma(t)}=\langle\cdot,\cdot\rangle_t,\ t\in(0,b).$
\end{proposition}
\begin{proof} Let $g$ be a $G$-invariant metric on $M_0.$ Then, for each $h\in H,$ $\phi_h^*g=g.$ Evaluating at $\gamma(t),$ we obtain
\[
g_{\gamma(t)}\left((d\phi_h)_{\gamma(t)}(u),(d\phi_h)_{\gamma(t)}(v)\right)=g_{\gamma(t)}(u,v),\ u,v\in T_{\gamma(t)}M.
\]
Equivalently, 
\[
g_{\gamma(t)}(\rho_t(h)u,\rho_t(h)v)=g_{\gamma(t)}(u,v),\ u,v\in T_{\gamma(t)}M.
\]
Thus, for each $t\in(0,b)$, $g_{\gamma(t)}$ is a $\rho_t(H)$-invariant inner product on $T_{\gamma(t)}M$. Given a basis $\{X_1,...,X_n\}$ of $\mathfrak{m},$ define
\[
\widetilde{E}_0:=\widetilde{\gamma'},\ \widetilde{E}_j:=X_j^*,\ j=1,...,n,
\]
and let $E_j(t):=\widetilde{E}_j(\gamma(t)),\ j=0,...,n.$ Since $g$ is a Riemannian metric, the functions $g(\widetilde{E}_j,\widetilde{E}_k),\ j,k=0,...,n,$ are smooth. Therefore, the functions
\[
g(\widetilde{E}_j,\widetilde{E}_k)\circ\gamma,\ j,k=0,...,n,
\]
are smooth. These are precisely the functions in \eqref{eq:metric-Ej-Ek}.\\

Conversely, let $(\langle\cdot,\cdot\rangle_t)_{t\in(0,b)}$ be a family such that, for each $t\in(0,b),$ $\langle \cdot,\cdot\rangle_t$ is a $\rho_t(H)$-invariant inner product on $T_{\gamma(t)}M$ and assume that $(\langle\cdot,\cdot\rangle_t)_{t\in(0,b)}$ is smooth. Given $x\in M_0$ and $u,v\in T_x M_0,$ there exists $a\in G$ such that $a\cdot \gamma(t(x))=x.$ Define 
\[
g_x(u,v):=\langle(d\phi_{a^{-1}})_{x}(u),(d\phi_{a^{-1}})_{x}(v)\rangle_{t(x)}.
\]
This definition does not depend on the choice of $a.$ In fact, if $a_1\cdot\gamma(t(x))=a_2\cdot\gamma(t(x)),$ then $a_2^{-1}a_1=h\in H.$ Since $\langle\cdot,\cdot\rangle_{t(x)}$ is $\rho_{t(x)}(H)$-invariant,
\begin{align*}
\langle(d\phi_{a_1^{-1}})_{x}(u),(d\phi_{a_1^{-1}})_{x}(v)\rangle_{t(x)}&=\langle\rho_{t(x)}(h)(d\phi_{a_1^{-1}})_{x}(u),\rho_{t(x)}(h)(d\phi_{a_1^{-1}})_{x}(v)\rangle_{t(x)}\\
&=\langle(d\phi_{h})_{\gamma(t(x))}(d\phi_{a_1^{-1}})_x(u),(d\phi_{h})_{\gamma(t(x))}(d\phi_{a_1^{-1}})_x(v)\rangle_{t(x)}\\
&=\langle(d\phi_{ha_1^{-1}})_{x}(u),(d\phi_{ha_1^{-1}})_{x}(v)\rangle_{t(x)}\\
&=\langle(d\phi_{a_2^{-1}})_{x}(u),(d\phi_{a_2^{-1}})_{x}(v)\rangle_{t(x)}.
\end{align*}
By construction, $g$ is a $G$-invariant field of inner products on $M_0$ and is the unique one satisfying $g_{\gamma(t)}=\langle\cdot,\cdot\rangle_t,\ t\in(0,b).$ To show that $g$ is smooth, let $\{X_1,...,X_n\}$ be a basis of $\mathfrak{m}$ and $E_0(t):=\gamma'(t),\ E_j(t):=X_j^*(\gamma(t)),\ t\in(0,b),\ j=1,...,n.$ Let $U\subseteq G/H$ be an open set admitting a smooth local section $\sigma:U\longrightarrow G$ of the projection $G\to G/H$. For $j=0,...,n$, define a vector field $\widehat{E}_j$ on $\Psi((0,b)\times U)$ by
\[
\widehat{E}_j(\Psi(t,xH)):=(d\phi_{\sigma(xH)})_{\gamma(t)}E_j(t),\ (t,xH)\in(0,b)\times U.
\]
The vector fields $\widehat{E}_0,...,\widehat{E}_n$ form a smooth local frame on $\Psi((0,b)\times U)$. By the definition of $g$,
\[
g_{\Psi(t,xH)}\left(\widehat{E}_j(\Psi(t,xH)),\widehat{E}_k(\Psi(t,xH))\right)=\langle E_j(t),E_k(t)\rangle_t.
\]
Since the family $(\langle\cdot,\cdot\rangle_t)_{t\in(0,b)}$ is smooth, the functions $t\longmapsto\langle E_j(t),E_k(t)\rangle_t$ are smooth. Hence, the functions 
\[
(t,xH)\longmapsto \langle E_j(t),E_k(t)\rangle_t
\]
are smooth on $(0,b)\times U$. Therefore, $g$ is smooth on $\Psi((0,b)\times U)$. Since $G/H$ can be covered by open sets of this form, $g$ is smooth on $M_0$.
\end{proof}
\begin{remark}
Given two bases $\{X_1,...,X_n\}$ and $\{Y_1,...,Y_n\}$ of $\mathfrak{m}$, there exists a matrix $C\in\mathrm{GL}(n,\mathbb{R})$ such that
\[
(Y_1\ \cdots\ Y_n)=(X_1\ \cdots\ X_n)\,C.
\]
Consequently, for each $t\in(0,b)$,
\[
\left(\gamma'(t)\ \ Y_1^*(\gamma(t))\ \ \cdots\ \ Y_n^*(\gamma(t))\right)=\left(\gamma'(t)\ \ X_1^*(\gamma(t))\ \ \cdots\ \ X_n^*(\gamma(t))\right)
\mathrm{diag}(1,C).
\]
Thus, the change-of-basis matrix between 
\[
\{\gamma'(t),X_1^*(\gamma(t)),...,X_n^*(\gamma(t))\}
\]
and 
\[
\{\gamma'(t),Y_1^*(\gamma(t)),...,Y_n^*(\gamma(t))\}
\]
does not depend on $t$. Therefore, a family $(\langle\cdot,\cdot\rangle_t)_{t\in(0,b)}$ is smooth if and only if, for one basis $\{X_1,...,X_n\}$ of $\mathfrak{m}$, the functions
\[
t\longmapsto
\langle E_j(t),E_k(t)\rangle_t,\ j,k=0,...,n,
\]
are smooth. \qed
\end{remark}
We consider the particular family  $((\cdot,\cdot)_t)_{t\in(0,b)}$ of inner products defined by
\begin{align}\label{eq:fixed-family-inner-products}
\begin{split}
&(\gamma'(t),\gamma'(t))_t:=1,\\
&(X^*(\gamma(t)),Y^*(\gamma(t)))_t:=(X,Y),\\
&(\gamma'(t),X^*(\gamma(t)))_t:=0,
\end{split}
\end{align}
for all $X,Y\in\mathfrak{m}.$ It is clear that this family is smooth and that each $(\cdot,\cdot)_t$ is $\rho_t(H)$-invariant. By the previous proposition, for any $G$-invariant metric $g$ on $M_0,$ $g_{\gamma(t)}$ is a $\rho_t(H)$-invariant inner product on $T_{\gamma(t)}M.$ Therefore, there exists a family $(\Lambda(t))_{t\in(0,b)}$ of $\rho_t$-equivariant, positive and $(\cdot,\cdot)_t$-symmetric operators $\Lambda(t):T_{\gamma(t)}M\longrightarrow T_{\gamma(t)}M$ such that
\[
g_{\gamma(t)}(u,v)=(\Lambda(t)u,v)_t,
\]
for all $u,v\in T_{\gamma(t)}M.$ The family $(\Lambda(t))_{t\in(0,b)}$ is smooth in the sense that, for one, and hence any, basis $\{X_1,...,X_n\}$ of $\mathfrak{m},$ setting $\mathcal{B}_t:=\{\gamma'(t),X_1^*(\gamma(t)),...,X_n^*(\gamma(t))\},$ the map
\[
t\longmapsto\left[\Lambda(t)\right]_{\mathcal{B}_t}\in M_{n+1}(\mathbb{R})
\]
is smooth. By Proposition \ref{prop:invariant-metrics}, any smooth family $(\Lambda(t))_{t\in(0,b)}$ of $\rho_t$-equivariant, positive and $(\cdot,\cdot)_t$-symmetric operators $\Lambda(t):T_{\gamma(t)}M\longrightarrow T_{\gamma(t)}M$ determines a unique $G$-invariant Riemannian metric on $M_0.$

\subsection{Smoothness at singular orbits} We now recall the setting used by Verdiani and Ziller \cite{VZ} to study the smooth extension of invariant metrics to a singular orbit. With the notations of the previous subsection, we have that the transitive action of $K_0$ on $K_0/H\approx \mathbb{S}^{\ell_0}\subseteq \mathbb{R}^{\ell_0+1}$ can be identified with the restriction to the unit sphere of a linear action of $K_0$ on $\mathbb R^{\ell_0+1}.$ By the slice theorem, there exists a $G$-invariant neighborhood $U_0$ of $Q_0$ in $M$ which is $G$-equivariantly diffeomorphic to $G\times_{K_0}D_{\delta},$ where $D_{\delta}\subseteq\mathbb{R}^{\ell_0+1}$ is a $K_0$-invariant disk and the zero section corresponds to $Q_0$. Thus, to study the smooth extension at $Q_0$, we may work on this homogeneous disk bundle. It is convenient to regard it as a neighborhood of the zero section in the homogeneous vector bundle
\[
E:=G\times_{K_0}\mathbb{R}^{\ell_0+1}.
\]
Here, we identify the zero section $G\times_{K_0}\{0\}$ with $Q_0$. The group $G$ acts on $E$ with cohomogeneity one by
\[
b\cdot[a,v]:=[ba,v].
\]
For $X\in\mathfrak g$, we denote by $X^\#$ the fundamental vector field on $E$ associated with $X$, namely,
\[
X^\#([a,v]):=\left.\frac{d}{dt}[\exp(tX)a,v]\,\right|_{t=0}.
\]
Let $q_0:=[e,0]\in Q_0.$ Its isotropy group is $G_{q_0}=K_0.$ For a nonzero vector $v\in \mathbb{R}^{\ell_0+1},$ whose isotropy group in $K_0$ is $H$, the point $[e,v]$ has isotropy group $H$. Hence, the principal orbits are diffeomorphic to $G/H$. Choose $v_0\in\mathbb{S}^{\ell_0}$ such that $(K_0)_{v_0}=H$ and consider the smooth curve $c:[0,\varepsilon)\longrightarrow E$ defined by $c(s):=[e,sv_0].$ Let $g$ be a smooth $G$-invariant Riemannian metric on $E\setminus Q_0$ for which $c|_{(0,\varepsilon)}$ is a unit-speed geodesic orthogonal to the principal orbits. For $s\in(0,\varepsilon)$, consider the isotropy representation  $\rho_s:H\longrightarrow \mathrm{GL}(T_{c(s)}E).$ As before, $g$ can be written as
\[
g=ds^2+\langle\cdot,\cdot\rangle_s,\ s\in(0,\varepsilon),
\]
where $(\langle\cdot,\cdot\rangle_s)_{s\in(0,\varepsilon)}$ is a smooth family such that $\langle\cdot,\cdot\rangle_s$ is a $\rho_s(H)$-invariant inner product on $T_{c(s)}(G\cdot c(s)).$\\

With respect to the fixed $\operatorname{Ad}(G)$-invariant inner product
$(\cdot,\cdot)$ on $\mathfrak{g}$, let $\mathfrak{n}_0$ be the
orthogonal complement of $\mathfrak{k}_0$ in $\mathfrak{g}$ and
$\mathfrak{p}_0$ the orthogonal complement of $\mathfrak{h}$ in
$\mathfrak{k}_0$. Thus,
\[
\mathfrak{g}=\mathfrak{k}_0\oplus\mathfrak{n}_0,\ \mathfrak{k}_0=\mathfrak{h}\oplus\mathfrak{p}_0,\ \mathfrak{m}=\mathfrak{p}_0\oplus\mathfrak{n}_0.
\]
Assume that the isotropy representation of $H$ on $\mathfrak{m}$ decomposes as
\[
\mathfrak{m}=\mathfrak{m}_1\oplus\cdots\oplus\mathfrak{m}_{\nu},
\]
where each $\mathfrak{m}_j$ is an irreducible $\operatorname{Ad}(H)$-module and
$\mathfrak{m}_j$ and $\mathfrak{m}_k$ are inequivalent whenever $j\neq k$.\footnote{The assumption that the modules are pairwise inequivalent is not made in \cite{VZ}. We impose it here because it holds for the action considered in this paper and simplifies the description of the invariant metrics and the smoothness conditions.} Since $\mathfrak{p}_0$ and $\mathfrak{n}_0$ are $\operatorname{Ad}(H)$-invariant, after reordering the irreducible summands, we may write
\[
\mathfrak{p}_0=\mathfrak{m}_1\oplus\cdots\oplus\mathfrak{m}_{\mu},\ \mathfrak{n}_0=
\mathfrak{m}_{\mu+1}\oplus\cdots\oplus\mathfrak{m}_{\nu}.
\]
We set $\mathfrak{p}_0^\beta:=\mathfrak{m}_\beta,\ \beta=1,...,\mu.$ Hence,
\[
\mathfrak{m}=\mathfrak{p}_0^1\oplus\cdots\oplus\mathfrak{p}_0^\mu\oplus\mathfrak{n}_0.
\]
For each $\beta=1,...,\mu$, choose $X_\beta\in\mathfrak{p}_0^\beta$ generating a closed one-parameter subgroup $L_\beta:=\{\exp(\theta X_\beta):\theta\in\mathbb{R}\}\subseteq K_0,$ normalized so that $\theta\mapsto\exp(\theta X_\beta)$ has period $2\pi,$ and set
\[
a_\beta:=|L_\beta\cap H|.
\]
For each $\beta$, decompose $\mathfrak{n}_0$ as
\[
\mathfrak{n}_0=\mathfrak{l}_{\beta,0}\oplus\mathfrak{l}_{\beta,1}\oplus\cdots\oplus\mathfrak{l}_{\beta,N_\beta},
\]
where $L_\beta$ acts trivially on $\mathfrak{l}_{\beta,0}$ and each $\mathfrak{l}_{\beta,\alpha}$, $\alpha=1,...,N_\beta$, is a two-dimensional irreducible real $L_\beta$-module. Fix a basis $\{U_{\beta,1},...,U_{\beta,r_\beta}\}$ of $\mathfrak{l}_{\beta,0}.$ For each $\alpha\in\{1,...,N_\beta\},$ choose a basis $\{Y_{1,\beta,\alpha},Y_{2,\beta,\alpha}\}$ of $\mathfrak{l}_{\beta,\alpha}$ such that
\[
[X_\beta,Y_{1,\beta,\alpha}]=d_{\beta,\alpha}Y_{2,\beta,\alpha},\ [X_\beta,Y_{2,\beta,\alpha}]=-d_{\beta,\alpha}Y_{1,\beta,\alpha},
\]
with $d_{\beta,\alpha}\in\mathbb{N}.$ The conditions of \cite{VZ} can then be written directly in terms of these decompositions. 
\begin{proposition}[\cite{VZ}]\label{prop:V-Z-smoothness-conditions}
The metric tensor $g$ extends to a smooth $G$-invariant symmetric $2$-tensor on a neighborhood of $Q_0$ if and only if for every $\beta=1,...,\mu,$
\[
\begin{aligned}
&g_{c(s)}(X_\beta^\#,X_\beta^\#)=a_\beta^2s^2+s^4\phi_\beta(s^2),\\[0.4em]
&g_{c(s)}(U_{\beta,j}^\#,U_{\beta,k}^\#)=\phi_{j,k,\beta}(s^2),\ 1\leq j\leq k\leq r_\beta,\\[0.4em]
&g_{c(s)}(U_{\beta,j}^\#,Y_{\ell,\beta,\alpha}^\#)=s^{d_{\beta,\alpha}/a_\beta}
\phi_{j,\ell,\beta,\alpha}(s^2),\ j=1,...,r_\beta,\ \ell=1,2,\ \alpha=1,...,N_\beta,\\[0.4em]
&g_{c(s)}(Y_{1,\beta,\alpha}^\#,Y_{1,\beta,\alpha}^\#)+g_{c(s)}(Y_{2,\beta,\alpha}^\#,Y_{2,\beta,\alpha}^\#)=\phi_{1,\beta,\alpha}(s^2),\ \alpha=1,...,N_\beta,\\[0.4em]
&g_{c(s)}(Y_{1,\beta,\alpha}^\#,Y_{1,\beta,\alpha}^\#)-g_{c(s)}(Y_{2,\beta,\alpha}^\#,Y_{2,\beta,\alpha}^\#)=s^{2d_{\beta,\alpha}/a_\beta}\phi_{2,\beta,\alpha}(s^2),\ \alpha=1,...,N_\beta,\\[0.4em]
&g_{c(s)}(Y_{1,\beta,\alpha}^\#,Y_{2,\beta,\alpha}^\#)=s^{2d_{\beta,\alpha}/a_\beta}\phi_{3,\beta,\alpha}(s^2),\ \alpha=1,...,N_\beta,\\[0.4em]
&g_{c(s)}(Y_{1,\beta,\alpha}^\#,Y_{1,\beta,\gamma}^\#)+g_{c(s)}(Y_{2,\beta,\alpha}^\#,Y_{2,\beta,\gamma}^\#)=s^{|d_{\beta,\alpha}-d_{\beta,\gamma}|/a_\beta}\phi_{1,\beta,\alpha,\gamma}(s^2),\\[0.4em]
&g_{c(s)}(Y_{1,\beta,\alpha}^\#,Y_{2,\beta,\gamma}^\#)-g_{c(s)}(Y_{2,\beta,\alpha}^\#,Y_{1,\beta,\gamma}^\#)=s^{|d_{\beta,\alpha}-d_{\beta,\gamma}|/a_\beta}\phi_{2,\beta,\alpha,\gamma}(s^2),\\[0.4em]
&g_{c(s)}(Y_{1,\beta,\alpha}^\#,Y_{1,\beta,\gamma}^\#)-g_{c(s)}(Y_{2,\beta,\alpha}^\#,Y_{2,\beta,\gamma}^\#)=s^{(d_{\beta,\alpha}+d_{\beta,\gamma})/a_\beta}\phi_{3,\beta,\alpha,\gamma}(s^2),\\[0.4em]
&g_{c(s)}(Y_{1,\beta,\alpha}^\#,Y_{2,\beta,\gamma}^\#)+g_{c(s)}(Y_{2,\beta,\alpha}^\#,Y_{1,\beta,\gamma}^\#)=s^{(d_{\beta,\alpha}+d_{\beta,\gamma})/a_\beta}\phi_{4,\beta,\alpha,\gamma}(s^2),\ 1\leq \alpha <\gamma\leq N_\beta,
\end{aligned}
\]
where the functions $\phi_\beta,\phi_{j,k,\beta},\phi_{j,\ell,\beta,\alpha},\phi_{j,\beta,\alpha},\phi_{j,\beta,\alpha,\gamma}$ are smooth on an interval of the form $[0,\delta)$ and, if any exponent appearing above is not an integer, the corresponding left-hand side must vanish identically. The smooth extension of $g$ is a Riemannian metric if and only if it is positive definite along $Q_0.$
\end{proposition}

\section{The $G_{p,q}$-action}\label{section:Gpq-action}
Fix integers $p\geq q\geq 1$ with $p>1.$ Let $G_{p,q}:=\mathrm{U}(p+1)\times\mathrm{U}(q+1)$ and
\[
\mathbb{C}P^{p+q+1}=\{[z:w]:z\in\mathbb{C}^{p+1},\ w\in\mathbb{C}^{q+1},\ (z,w)\neq(0,0)\}.
\]
Consider the cohomogeneity one action
\[
\begin{aligned}
\phi:G_{p,q}\times\mathbb{C}P^{p+q+1}&\longrightarrow\mathbb{C}P^{p+q+1},\\
((U,V),[z:w])&\longmapsto[Uz:Vw].
\end{aligned}
\]
Let $\{e_1,...,e_{p+1}\}$ and $\{f_1,...,f_{q+1}\}$ denote the canonical bases of $\mathbb{C}^{p+1}$ and $\mathbb{C}^{q+1}$, respectively. The curve
\[
\gamma(t):=\big[(\cos t)e_{p+1}:(\sin t)f_{q+1}\big],\ 0\leq t\leq\frac{\pi}{2},
\]
intersects every $G_{p,q}$-orbit in a single point. For $0<t<\frac{\pi}{2}$, the isotropy subgroup of $\gamma(t)$ is
\[
H=\left\{\left(\mathrm{diag}(U_1,\lambda),\mathrm{diag}(V_1,\lambda)\right):U_1\in\mathrm{U}(p),\ V_1\in\mathrm{U}(q),\ \lambda\in\mathrm{U}(1)\right\}\cong\mathrm{U}(p)\times\mathrm{U}(q)\times\mathrm{U}(1).
\]
Hence, the principal orbits are diffeomorphic to $G_{p,q}/H.$ The isotropy groups of $\gamma(0)$ and $\gamma(\pi/2)$ are, respectively,
\[
K_0=\left\{(U,V)\in G_{p,q}:U=\mathrm{diag}(U_1,\lambda),\ U_1\in\mathrm{U}(p),\ \lambda\in\mathrm{U}(1)\right\}\cong\mathrm{U}(p)\times\mathrm{U}(1)\times\mathrm{U}(q+1)
\]
and
\[
K_1=\left\{(U,V)\in G_{p,q}:V=\mathrm{diag}(V_1,\lambda),\ V_1\in\mathrm{U}(q),\ \lambda\in\mathrm{U}(1) \right\}\cong\mathrm{U}(p+1)\times\mathrm{U}(q)\times\mathrm{U}(1).
\]
Therefore, the singular orbits corresponding to $t=0$ and $t=\pi/2$ are, respectively,
\[
Q_0=G_{p,q}/K_0\approx\frac{\mathrm{U}(p+1)}{\mathrm{U}(p)\times\mathrm{U}(1)}\approx\mathbb{C}P^p
\]
and
\[
Q_1=G_{p,q}/K_1\approx\frac{\mathrm{U}(q+1)}{\mathrm{U}(q)\times\mathrm{U}(1)}\approx\mathbb{C}P^q.
\]

Let $\mathfrak{g}_{p,q}:=\mathfrak{u}(p+1)\oplus\mathfrak{u}(q+1)$ and endow $\mathfrak{g}_{p,q}$ with the $\operatorname{Ad}(G_{p,q})$-invariant inner product $(\cdot,\cdot)$ defined by 
\begin{equation*}
\left((Z_1,Z_2),(Z_3,Z_4)\right)=\frac{1}{2}\left[\mathrm{Tr}(Y_1Y_3-X_1X_3)+\mathrm{Tr}(Y_2Y_4-X_2X_4)\right],
\end{equation*}
where $Z_j=X_j+iY_j,\ j=1,2,3,4$ and $X_j,Y_j$ have real entries. Let $\mathfrak{h}\subseteq \mathfrak{g}_{p,q}$ be the Lie algebra of $H$ and let $\mathfrak{m}$ be the $(\cdot,\cdot)$-orthogonal complement of $\mathfrak{h}$ in $\mathfrak{g}_{p,q}.$ Then $\mathfrak{m}$ decomposes into the three irreducible, pairwise inequivalent $H$-submodules
\begin{align*}
\mathfrak{m}_1&=\mathrm{span}_{\mathbb{R}}\left\{\left(iE_{p+1,p+1},-iF_{q+1,q+1}\right)\right\},\\
\mathfrak{m}_2&=\mathrm{span}_{\mathbb{R}}\left\{
\left(E_{p+1,j}-E_{j,p+1},0\right),\left(i(E_{p+1,j}+E_{j,p+1}),0\right):1\leq j\leq p\right\},\\
\mathfrak{m}_3&=\mathrm{span}_{\mathbb{R}}\left\{\left(0,F_{q+1,j}-F_{j,q+1}\right),\left(0,i(F_{q+1,j}+F_{j,q+1})\right):1\leq j\leq q\right\},
\end{align*}
where $E_{a,b}$ denotes the $(p+1)\times(p+1)$-matrix with a one in its $(a,b)$-entry and zero elsewhere, and $F_{c,d}$ denotes the $(q+1)\times(q+1)$-matrix with a one in its $(c,d)$-entry and zero elsewhere. Define
\begin{align*}
    &Z:=(iE_{p+1,p+1},-iF_{q+1,q+1}),\\[0.2em]
    &A_j:=(E_{p+1,j}-E_{j,p+1},0),\ B_j:=(i(E_{p+1,j}+E_{j,p+1}),0),\ j=1,...,p,\\[0.2em]
    &C_j:=(0,F_{q+1,j}-F_{j,q+1}),\ D_j:=(0,i(F_{q+1,j}+F_{j,q+1})),\ j=1,...,q. 
\end{align*}
Then $\{Z\}\cup\{A_j,B_j:j=1,...,p\}\cup\{C_j,D_j:j=1,...,q\}$ is an $(\cdot,\cdot)$-orthonormal basis of $\mathfrak{m}$ adapted to the decomposition $\mathfrak{m}=\mathfrak{m}_1\oplus\mathfrak{m}_2\oplus\mathfrak{m}_3.$ For each $t\in(0,\pi/2),$ under the identification $\mathfrak{m}\simeq T_{\gamma(t)}P_t$ introduced in Subsection~\ref{subsec:invariant-metrics}, we have the decomposition
\[
T_{\gamma(t)}\mathbb{C}P^{p+q+1}=\mathbb{R}\gamma'(t)\oplus\mathfrak{m}_1\oplus\mathfrak{m}_2\oplus\mathfrak{m}_3
\]
into irreducible $\rho_t(H)$-invariant subspaces. The $H$-modules $\mathbb{R}\gamma'(t)$ and $\mathfrak{m}_1$ are equivalent, with trivial $H$-action. On the other hand, the modules $\mathfrak{m}_2$ and $\mathfrak{m}_3$ are nontrivial and inequivalent. Consequently, for each $t\in(0,\pi/2),$ every positive, $(\cdot,\cdot)_t$-symmetric and $\rho_t$-equivariant operator $\Lambda(t):T_{\gamma(t)}\mathbb{C}P^{p+q+1}\longrightarrow T_{\gamma(t)}\mathbb{C}P^{p+q+1}$ is of the form 
\begin{align*}
    &\Lambda(t)\gamma'(t)=f_0(t)^2\gamma'(t)+u(t)Z^*(\gamma(t)),\\[0.2em]
    &\Lambda(t)Z^*(\gamma(t))=u(t)\gamma'(t)+f_1(t)^2Z^*(\gamma(t)),\\[0.2em]
    &\Lambda(t)A_j^*(\gamma(t))=f_2(t)^2A_j^*(\gamma(t)),\ \Lambda(t)B_j^*(\gamma(t))=f_2(t)^2B_j^*(\gamma(t)),\ j=1,...,p,\\[0.2em]
    &\Lambda(t)C_j^*(\gamma(t))=f_3(t)^2C_j^*(\gamma(t)),\ \Lambda(t)D_j^*(\gamma(t))=f_3(t)^2D_j^*(\gamma(t)),\ j=1,...,q,
\end{align*}
with $f_k(t)>0,\ k=0,1,2,3$ and $|u(t)|< f_0(t)f_1(t).$  By the discussion of Subsection~\ref{subsec:invariant-metrics}, we have that the set of all $G_{p,q}$-invariant metrics $g$ on $(\mathbb{C}P^{p+q+1})_0:=\mathbb{C}P^{p+q+1}\setminus(Q_0\cup Q_1)$ is in bijection with the set of smooth maps
\[
\left(0,\frac{\pi}{2}\right)\ni t\longmapsto (f_0(t),f_1(t),f_2(t),f_3(t),u(t))
\]
such that $f_k(t)>0,\ k=0,1,2,3,$ and $|u(t)|<f_0(t)f_1(t),$ for all $t\in(0,\pi/2).$ The correspondence is given by the formulas
\begin{align*}
    &g_{\gamma(t)}\left(\gamma'(t),\gamma'(t)\right)=f_0(t)^2,\\[0.2em]
    &g_{\gamma(t)}\left(Z^*(\gamma(t)),Z^*(\gamma(t))\right)=f_1(t)^2,\\[0.2em]
    &g_{\gamma(t)}\left(A_j^*(\gamma(t)),A_j^*(\gamma(t))\right)=g_{\gamma(t)}\left(B_j^*(\gamma(t)),B_j^*(\gamma(t))\right)=f_2(t)^2,\ j=1,...,p,\\[0.2em]
    &g_{\gamma(t)}\left(C_j^*(\gamma(t)),C_j^*(\gamma(t))\right)=g_{\gamma(t)}\left(D_j^*(\gamma(t)),D_j^*(\gamma(t))\right)=f_3(t)^2,\ j=1,...,q,\\[0.2em]
    &g_{\gamma(t)}\left(\gamma'(t),Z^*(\gamma(t))\right)=g_{\gamma(t)}\left(Z^*(\gamma(t)),\gamma'(t)\right)=u(t),
\end{align*}
and all the remaining pairings between these vectors vanish. We will write 
\[
g=(f_0,f_1,f_2,f_3,u).
\]
Now, let $J$ be the canonical complex structure on $\mathbb{C}P^{p+q+1}$ and let
\[
\pi:\mathbb{C}^{p+q+2}\setminus\{0\}\longrightarrow\mathbb{C}P^{p+q+1},\ \pi(z):=[z]
\]
be the canonical projection. For each $z\in\mathbb{C}^{p+q+2}\setminus\{0\},$ we have $\ker(d\pi)_z=\mathbb{C}z.$ The complex structure $J$ is characterized by
\begin{equation}\label{eq:standard-complex-structure}
    J_{\pi(z)}((d\pi)_z(v))=(d\pi)_z(iv),\ z\in\mathbb{C}^{p+q+2}\setminus\{0\},\ v\in\mathbb{C}^{p+q+2}.
\end{equation}
Since the action of $G_{p,q}$ on $\mathbb{C}P^{p+q+1}$ is induced by complex linear transformations of $\mathbb{C}^{p+q+2}$, the complex
structure $J$ is $G_{p,q}$-invariant.

\begin{proposition}\label{prop:canonical-J}
For each $t\in(0,\pi/2)$, the canonical complex structure $J$ satisfies
\[
\begin{aligned}
&J\gamma'(t)=-\frac{1}{\sin(2t)}Z^*(\gamma(t)),\ JZ^*(\gamma(t))=\sin(2t)\gamma'(t),\\[0.3em]
&JA_j^*(\gamma(t))=-B_j^*(\gamma(t)),\ JB_j^*(\gamma(t))=A_j^*(\gamma(t)),\ j=1,...,p,\\[0.3em]
&JC_j^*(\gamma(t))=-D_j^*(\gamma(t)),\ JD_j^*(\gamma(t))=C_j^*(\gamma(t)),\ j=1,...,q.
\end{aligned}
\]
\end{proposition}
\begin{proof}
Consider the curve $z(t):=((\cos t)e_{p+1},(\sin t)f_{q+1})\in\mathbb{C}^{p+q+2}.$ Then 
\[
z'(t)=(-(\sin t)e_{p+1},(\cos t)f_{q+1})\ \text{and}\ \gamma'(t)=(d\pi)_{z(t)}(z'(t)).
\]
By \eqref{eq:standard-complex-structure}, we have
\[
J\gamma'(t)=(d\pi)_{z(t)}\left(iz'(t)\right).
\]

On the other hand, 
\begin{align*}
    Z^*(\gamma(t))&=\left.\frac{d}{ds}\exp(sZ)\cdot\gamma(t)\,\right|_{s=0}\\[0.2em]
    &=\left.\frac{d}{ds}\left(\exp(siE_{p+1,p+1}),\exp(-siF_{q+1,q+1})\right)\cdot[(\cos t)e_{p+1}:(\sin t)f_{q+1}]\,\right|_{s=0}\\[0.2em]
    &=\left.\frac{d}{ds}\left[(\cos t)\exp(siE_{p+1,p+1})e_{p+1}:(\sin t)\exp(-siF_{q+1,q+1})f_{q+1}\right]\,\right|_{s=0}\\[0.2em]
    &=\left.\frac{d}{ds}\pi\left((\cos t)\exp(siE_{p+1,p+1})e_{p+1},(\sin t)\exp(-siF_{q+1,q+1})f_{q+1}\right)\,\right|_{s=0}\\[0.2em]
    &=(d\pi)_{z(t)}\left(\left.\frac{d}{ds}\left((\cos t)\exp(siE_{p+1,p+1})e_{p+1},(\sin t)\exp(-siF_{q+1,q+1})f_{q+1}\right)\,\right|_{s=0}\right)\\[0.2em]
    &=(d\pi)_{z(t)}\left(i\left((\cos t)e_{p+1},-(\sin t)f_{q+1}\right)\right).
\end{align*}
From the identity
\[
i\left((\cos t)e_{p+1},-(\sin t)f_{q+1}\right)=-\sin(2t)iz'(t)+\cos(2t)iz(t),
\]
and since $iz(t)\in\ker(d\pi)_{z(t)}=\mathbb{C}z(t)$, we obtain
\[
Z^*(\gamma(t))=-\sin(2t)(d\pi)_{z(t)}(iz'(t))=-\sin(2t)J\gamma'(t).
\]
Since $J^2=-\operatorname{Id}$, it follows that $JZ^*(\gamma(t))=\sin(2t)\gamma'(t).$
Therefore,
\[
J\gamma'(t)=-\frac{1}{\sin(2t)}Z^*(\gamma(t)),\ JZ^*(\gamma(t))=\sin(2t)\gamma'(t).
\]
Now observe that, for $j=1,...,p,$
\begin{align*}
    A_j^*(\gamma(t))&=(d\pi)_{z(t)}\left(\left.\frac{d}{ds}\left((\cos t)\exp(s(E_{p+1,j}-E_{j,p+1}))e_{p+1},(\sin t)f_{q+1}\right)\,\right|_{s=0}\right)\\[0.2em]
    &=(d\pi)_{z(t)}(-(\cos t)e_j,0),\\[0.2em]
    B_j^*(\gamma(t))&=(d\pi)_{z(t)}\left(\left.\frac{d}{ds}\left((\cos t)\exp(si(E_{p+1,j}+E_{j,p+1}))e_{p+1},(\sin t)f_{q+1}\right)\,\right|_{s=0}\right)\\[0.2em]
    &=(d\pi)_{z(t)}(i(\cos t)e_j,0).
\end{align*}
Therefore
\[
JA_j^*(\gamma(t))=-B_j^*(\gamma(t)),\ JB_j^*(\gamma(t))=A_j^*(\gamma(t)),\ j=1,...,p.
\]
Analogously,
\begin{align*}
    C_j^*(\gamma(t))=(d\pi)_{z(t)}(0,-(\sin t)f_j),\ D_j^*(\gamma(t))=(d\pi)_{z(t)}(0,i(\sin t)f_j),\ j=1,...,q,
\end{align*}
thus
\[
JC_j^*(\gamma(t))=-D_j^*(\gamma(t)),\ JD_j^*(\gamma(t))=C_j^*(\gamma(t)),\ j=1,...,q.
\]
\end{proof}
\begin{corollary}\label{coro:compatible-metrics}
Let $g=(f_0,f_1,f_2,f_3,u)$ be a $G_{p,q}$-invariant Riemannian metric on $(\mathbb{C}P^{p+q+1})_0.$ Then $g$ is compatible with $J$ if and only if
\[
u(t)=0,\ f_1(t)=\sin(2t)f_0(t),\ t\in(0,\pi/2).
\]
\end{corollary}
\begin{proof}
Assume that $g$ and $J$ are compatible. Then, for each $t\in(0,\pi/2),$ we have
\begin{align*}
    u(t)&=g_{\gamma(t)}\left(\gamma'(t),Z^*(\gamma(t))\right)\\[0.2em]
        &=g_{\gamma(t)}\left(J\gamma'(t),JZ^*(\gamma(t))\right)\\[0.2em]
        &=g_{\gamma(t)}\left(-\frac{1}{\sin(2t)}Z^*(\gamma(t)),\sin(2t)\gamma'(t)\right)\\[0.2em]
        &=-g_{\gamma(t)}\left(Z^*(\gamma(t)),\gamma'(t)\right)\\[0.2em]
        &=-u(t).
\end{align*}    
Thus, $u(t)=0.$ Also,
\begin{align*}
    f_1(t)^2&=g_{\gamma(t)}\left(Z^*(\gamma(t)),Z^*(\gamma(t))\right)\\[0.2em]
            &=g_{\gamma(t)}\left(JZ^*(\gamma(t)),JZ^*(\gamma(t))\right)\\[0.2em]
            &=g_{\gamma(t)}\left(\sin(2t)\gamma'(t),\sin(2t)\gamma'(t)\right)\\[0.2em]
            &=\sin^2(2t)f_0(t)^2.
\end{align*}
Therefore, $f_1(t)=\sin(2t)f_0(t),\ t\in(0,\pi/2).$\\

Conversely, assume that $u(t)=0$ and $f_1(t)=\sin(2t)f_0(t),\ t\in(0,\pi/2).$ Since both $g$ and $J$ are $G_{p,q}$-invariant, it is enough to verify compatibility at the points $\gamma(t),\ t\in(0,\pi/2).$ A direct
computation using Proposition~\ref{prop:canonical-J} shows that $g_{\gamma(t)}(Jv_1,Jv_2)=g_{\gamma(t)}(v_1,v_2)$ for all $v_1,v_2$ in the basis
\[
\mathcal{B}_t=\left\{\gamma'(t),Z^*(\gamma(t)),A_j^*(\gamma(t)),B_j^*(\gamma(t)),C_k^*(\gamma(t)),D_k^*(\gamma(t)):j=1,...,p,\ k=1,...,q\right\}
\]
of $T_{\gamma(t)}\mathbb{C}P^{p+q+1}$. By bilinearity, $g_{\gamma(t)}(Ju_1,Ju_2)=g_{\gamma(t)}(u_1,u_2)$ for all $u_1,u_2\in T_{\gamma(t)}\mathbb{C}P^{p+q+1}$. Hence, $g$ is compatible with $J$.
\end{proof}
Set
\[
\mathfrak{p}_0:=\mathfrak{m}_1\oplus\mathfrak{m}_3,\ \mathfrak{n}_0:=\mathfrak{m}_2,\ \mathfrak{p}_1:=\mathfrak{m}_1\oplus\mathfrak{m}_2,\ \mathfrak{n}_1:=\mathfrak{m}_3.
\]
If $\mathfrak{k}_0$ and $\mathfrak{k}_1$ denote the Lie algebras of $K_0$ and $K_1$, respectively, then
\[
\mathfrak{k}_0=\mathfrak{h}\oplus\mathfrak{p}_0,\ \mathfrak{k}_1=\mathfrak{h}\oplus\mathfrak{p}_1,\ \mathfrak{g}_{p,q}=\mathfrak{k}_0\oplus\mathfrak{n}_0,\ \mathfrak{g}_{p,q}=\mathfrak{k}_1\oplus\mathfrak{n}_1.
\]
These are $(\cdot,\cdot)$-orthogonal reductive decompositions, and
\[
K_0/H\approx\mathbb{S}^{2q+1},\ K_1/H\approx\mathbb{S}^{2p+1}.
\]
We now study when a $G_{p,q}$-invariant metric 
\[
g=(f_0,\sin(2t)f_0,f_2,f_3,0)
\]
defined on $(\mathbb{C}P^{p+q+1})_0$ extends to a smooth $G_{p,q}$-invariant Riemannian metric on $\mathbb{C}P^{p+q+1}.$ First observe that, if $g$ extends to a smooth Riemannian metric on a neighborhood of $Q_0$, then the function $t\longmapsto g_{\gamma(t)}(\gamma'(t),\gamma'(t))$ extends to a smooth even function near $t=0$ because
\[
\left(\mathrm{Id}_{p+1},\mathrm{diag}(\mathrm{Id}_q,-1)\right)\cdot\gamma(t)=\gamma(-t).
\]
Since $f_0(t)^2=g_{\gamma(t)}(\gamma'(t),\gamma'(t))$ for $t>0$ and the extended metric is positive definite, this smooth even extension is positive near $t=0.$ Its square root therefore gives a smooth even extension of $f_0.$ Hence, by \cite[Theorem~1]{Whitney1943}, $f_0(t)=\phi_0(t^2)$ for some smooth positive function $\phi_0$ near $0.$ The same argument at $Q_1$, using the variable $r=\frac{\pi}{2}-t$, yields $f_0\left(\frac{\pi}{2}-r\right)=\phi_1(r^2)$ for some smooth positive function $\phi_1$ near $0$. Consequently, from now on, we assume the following condition on $f_0$: there exist $\delta>0$ and smooth positive functions
$\phi_0,\phi_1:[0,\delta^2)\longrightarrow\mathbb{R}$ such that
\begin{equation}\label{eq:f0-smoothness}
f_0(t)=\phi_0(t^2),\ f_0\left(\frac{\pi}{2}-t\right)=\phi_1(t^2),\ 0<t<\delta.
\end{equation}
In particular, by defining $f_0(0):=\phi_0(0),$ $f_0\left(\frac{\pi}{2}\right):=\phi_1(0),$ the function $f_0$ extends to a smooth positive function on $[0,\pi/2].$
\begin{lemma}\label{lem:arc-length-Q0}
Define
\[
s_0(t):=\int_0^t f_0(r)\,dr.
\]
The following assertions hold:
\begin{itemize}
    \item[$a)$] For $\varepsilon>0$ sufficiently small, $s_0$ is a smooth diffeomorphism from $[0,\varepsilon)$ onto $[0,\varepsilon_0)$ for some $\varepsilon_0>0$. Let $\tau_0:=s_0^{-1}$ and set
    \[
        c_0(s):=\gamma(\tau_0(s)).
    \]
    Then $c_0$, with respect to $g$, is a unit-speed geodesic in $(\mathbb{C}P^{p+q+1})_0$ orthogonal to every principal orbit it meets.\\
    
    \item[$b)$] There exists a $G_{p,q}$-equivariant smooth diffeomorphism $\Phi_0$ from a neighborhood of the zero section in $G_{p,q}\times_{K_0}\mathbb{C}^{q+1}$ onto a neighborhood of $Q_0$ in $\mathbb{C}P^{p+q+1}$ such that
    \[
    \Phi_0\left(\left[(\mathrm{Id}_{p+1},\mathrm{Id}_{q+1}),sf_{q+1}\right]\right)=c_0(s),
    \]
    for all $s\geq 0$ sufficiently small.
\end{itemize}
\end{lemma}
\begin{proof} 
Since $f_0(t)=\phi_0(t^2)$ for $0<t<\delta$, where $\phi_0:[0,\delta^2)\longrightarrow\mathbb{R}$ is smooth and positive, the function $\widehat{f}_0(t):=\phi_0(t^2),\ -\delta<t<\delta,$ is a smooth positive even extension of $f_0$. Define
\[
\widehat{s}_0(t):=\int_0^t\widehat{f}_0(r)\,dr.
\]
Then $\widehat{s}_0$ is smooth and odd, and $\widehat{s}_0'(0)=\widehat{f}_0(0)=\phi_0(0)>0.$ Hence, by the inverse function theorem, there exist $\varepsilon,\varepsilon_0>0$ such that $\widehat{s}_0:(-\varepsilon,\varepsilon)\longrightarrow(-\varepsilon_0,\varepsilon_0)$ is a smooth diffeomorphism. Its restriction
\[
s_0:[0,\varepsilon)\longrightarrow[0,\varepsilon_0)
\]
is therefore a smooth diffeomorphism. Let $\widehat{\tau}_0:=\widehat{s}_0^{-1}$ and $\tau_0:=s_0^{-1}.$ Then $\tau_0$ is the restriction of $\widehat{\tau}_0$ to $[0,\varepsilon_0).$ Since $\widehat{s}_0$ is odd, $\widehat{\tau}_0$ is also odd. Moreover, $\widehat{\tau}_0(0)=0$, so we can write
\[
\widehat{\tau}_0(s)=s\eta_0(s),\ \eta_0(s):=\int_0^1\widehat{\tau}_0'(rs)\,dr.
\]
The function $\eta_0$ is smooth and even, so there exists a smooth function $\beta_0:[0,\varepsilon_0^2)\longrightarrow\mathbb{R}$ such that $\eta_0(s)=\beta_0(s^2).$ Since $\beta_0(0)=1/f_0(0)>0$, decreasing $\varepsilon_0$ if necessary, we may assume that $\beta_0$ is positive on $[0,\varepsilon_0^2).$\\

Now, let $c_0:[0,\varepsilon_0)\longrightarrow \mathbb{C}P^{p+q+1}$ be defined by $c_0(s):=\gamma(\tau_0(s)).$ For $s\in(0,\varepsilon_0)$, we have 
\[
c_0'(s)=\tau_0'(s)\gamma'(\tau_0(s))=\frac{1}{f_0(\tau_0(s))}\gamma'(\tau_0(s)).
\]
Therefore,
\begin{align*}
g_{c_0(s)}(c_0'(s),c_0'(s))=\frac{1}{f_0(\tau_0(s))^2}g_{\gamma(\tau_0(s))}\left(\gamma'(\tau_0(s)),\gamma'(\tau_0(s))\right)=1.
\end{align*}
Thus, $c_0$ has unit speed with respect to $g.$ Since the vector $\gamma'(t)$ is orthogonal to the tangent space of the principal orbit through $\gamma(t),$ $c_0'(s)$ is orthogonal to the principal orbit through $c_0(s)$ for every $s\in(0,\varepsilon_0).$ We now show that $c_0\big{|}_{(0,\varepsilon_0)}$ is geodesic. For each $X\in\mathfrak{m},$ the fundamental vector field $X^*$ is a Killing vector field because $g$ is $G_{p,q}$-invariant. Moreover,
\[
g_{c_0(s)}\left(c_0'(s),X^*(c_0(s))\right)=0,\ s\in(0,\varepsilon_0).
\]
Differentiating this identity with respect to $s$, we obtain
\[
g_{c_0(s)}\left(\left(\nabla_{c_0'}c_0'\right)(c_0(s)),X^*(c_0(s))\right)+g_{c_0(s)}\left(c_0'(s),\left(\nabla_{c_0'}X^*\right)(c_0(s))\right)=0.
\]
Since $X^*$ is Killing, $g_{c_0(s)}\left(c_0'(s),\left(\nabla_{c_0'}X^*\right)(c_0(s))\right)=0.$ Hence
\[
g_{c_0(s)}\left(\left(\nabla_{c_0'}c_0'\right)(c_0(s)),X^*(c_0(s))\right)=0.
\]
Thus $\left(\nabla_{c_0'}c_0'\right)(c_0(s))$ is orthogonal to the principal orbit through $c_0(s)$. On the other hand, $c_0$ has unit speed, thus
\[
g_{c_0(s)}\left(\left(\nabla_{c_0'}c_0'\right)(c_0(s)),c_0'(s)\right)=\frac{1}{2}\frac{d}{ds}
g_{c_0(s)}\left(c_0'(s),c_0'(s)\right)=0.
\]
Since $T_{c_0(s)}\mathbb{C}P^{p+q+1}=\mathbb{R}c_0'(s)\oplus T_{c_0(s)}P_{\tau_0(s)},$ we conclude that $\nabla_{c_0'}c_0'=0.$ Therefore, $c_0$ is a unit-speed geodesic in $(\mathbb{C}P^{p+q+1})_0$ orthogonal to every principal orbit it meets. This proves $a).$\\

Now, consider the unitary representation
\[
\sigma_0:K_0\longrightarrow \mathrm{U}(q+1),\ \sigma_0\left(\mathrm{diag}(U_1,\lambda),V\right):=\lambda^{-1}V,
\]
where $U_1\in\mathrm{U}(p)$, $\lambda\in\mathrm{U}(1)$ and $V\in\mathrm{U}(q+1)$. The manifold $G_{p,q}\times_{K_0}\mathbb{C}^{q+1}$ is the quotient of $G_{p,q}\times\mathbb{C}^{q+1}$ by the relation
\[
(a,w)\sim\left(ak,\sigma_0(k)^{-1}w\right),\ a\in G_{p,q},\ k\in K_0,\ w\in\mathbb{C}^{q+1}.
\]
Let $\mathcal{U}_0:=\mathbb{C}P^{p+q+1}\setminus Q_1=\left\{[z:w]\in\mathbb{C}P^{p+q+1}:z\neq0\right\}$ and define
\[
\Psi_0:G_{p,q}\times_{K_0}\mathbb{C}^{q+1}\ni\left[(U,V),w\right]\longmapsto[Ue_{p+1}:Vw]\in\mathcal{U}_0.
\]
Let us show that $\Psi_0$ is well defined. In fact, if $k=\left(\mathrm{diag}(U_1,\lambda),W\right)\in K_0,$ then $\sigma_0(k)^{-1}w=\lambda W^{-1}w.$ Consequently,
\begin{align*}
\Psi_0\left([(U,V)k,\sigma_0(k)^{-1}w]\right)&=\left[U\mathrm{diag}(U_1,\lambda)e_{p+1}:VW(\lambda W^{-1}w)\right]\\
&=[\lambda Ue_{p+1}:\lambda Vw]\\
&=[Ue_{p+1}:Vw]\\
&=\Psi_0\left([(U,V),w]\right).
\end{align*}
It is also $G_{p,q}$-equivariant because
\begin{align*}
    \Psi_0\left((\tilde{U},\tilde{V})\cdot[(U,V),w]\right)&=\Psi_0\left(\left[(\tilde{U}U,\tilde{V}V),w\right]\right)\\
    &=[\tilde{U}Ue_{p+1}:\tilde{V}Vw]\\
    &=(\tilde{U},\tilde{V})\cdot[Ue_{p+1}:Vw]\\
    &=(\tilde{U},\tilde{V})\cdot\Psi_0\left([(U,V),w]\right).
\end{align*}

The map $\Psi_0$ is smooth since it is induced by the smooth map
\[
G_{p,q}\times\mathbb{C}^{q+1}\ni((U,V),w)\longmapsto[Ue_{p+1}:Vw]\in\mathcal{U}_0.
\]
We now show that $\Psi_0$ is bijective. Given $[z:w]\in\mathcal{U}_0,$ choose $U_z\in\mathrm{U}(p+1)$ such that $U_ze_{p+1}=z/||z||,$ and define
\[
\Theta_0([z:w]):=\left[\left(U_z,\mathrm{Id}_{q+1}\right),\frac{w}{||z||}\right].
\]
This definition does not depend on the choice of $U_z$. In fact, if $U_z'e_{p+1}=z/||z||,$ then $U_z'=U_z\mathrm{diag}(U_1,1)$ for some $U_1\in\mathrm{U}(p)$ and $\left(\mathrm{diag}(U_1,1),\mathrm{Id}_{q+1}\right)\in K_0.$ Therefore
\begin{align*}
    \left[\left(U_z,\mathrm{Id}_{q+1}\right),\frac{w}{||z||}\right]=\left[(U_z\mathrm{diag}(U_1,1),\mathrm{Id}_{q+1}),1\mathrm{Id}_{q+1}^{-1}\left(\frac{w}{||z||}\right)\right]=\left[(U'_z,\mathrm{Id}_{q+1}),\frac{w}{||z||}\right].
\end{align*}
The definition is also independent of the representative of $[z:w].$ In fact, if $(z,w)$ is replaced by $(\mu z,\mu w),$ with $\mu\in\mathbb{C}\setminus\{0\}$ and $\lambda:=\mu/|\mu|\in\mathrm{U}(1),$ then
\[
\frac{\mu z}{||\mu z||}=\lambda\frac{z}{||z||},\ \frac{\mu w}{||\mu z||}=\lambda\frac{w}{||z||}.
\]
Taking $U_{\mu z}:=U_z\mathrm{diag}(\mathrm{Id}_p,\lambda)$ and using $\sigma_0\left(\mathrm{diag}(\mathrm{Id}_p,\lambda),\mathrm{Id}_{q+1}\right)^{-1}=\lambda\mathrm{Id}_{q+1},$ we obtain 
\[
\Theta_0([z:w])=\Theta_0([\mu z:\mu w]).
\]
Thus $\Theta_0$ is well defined. Clearly, $\Psi_0\circ\Theta_0=\mathrm{Id}_{\mathcal{U}_0}.$ On the other hand, if $\left[(U,V),w\right]\in G_{p,q}\times_{K_0}\mathbb{C}^{q+1},$ then
\[
\Theta_0\left(\Psi_0\left([(U,V),w]\right)\right)=\left[(U,\mathrm{Id}_{q+1}),Vw\right]=\left[(U,V),w\right],
\]
where the last equality follows by taking $k=\left(\mathrm{Id}_{p+1},V\right)\in K_0.$ Therefore, $\Theta_0=\Psi_0^{-1}$ and $\Psi_0$ is bijective.\\

Now, we show that $\Psi_0$ is a local diffeomorphism. By $G_{p,q}$-equivariance, it is enough to compute its differential at points of the form $\left[(\mathrm{Id}_{p+1},\mathrm{Id}_{q+1}),w\right].$ At such a point, the tangent space can be identified with the quotient $\left(\mathfrak{g}_{p,q}\oplus\mathbb{C}^{q+1}\right)/\sim,$ where
\[
((X,Y),\xi)\sim\left((X,Y)+A,\xi-(d\sigma_0)_{(\mathrm{Id}_{p+1},\mathrm{Id}_{q+1})}(A)w\right),\ A\in\mathfrak{k}_0.
\]

For $((X,Y),\xi)\in\mathfrak{g}_{p,q}\oplus\mathbb{C}^{q+1}$, consider the curve $t\longmapsto
\left[\left(\exp(tX),\exp(tY)\right),w+t\xi\right]$ through $\left[(\operatorname{Id}_{p+1},\operatorname{Id}_{q+1}),w\right]$. Then
\begin{align*}
(d\Psi_0)_{\left[(\mathrm{Id}_{p+1},\mathrm{Id}_{q+1}),w\right]}\left[(X,Y),\xi\right]&=\left.\frac{d}{dt}\Psi_0\left(\left[\left(\exp(tX),\exp(tY)\right),w+t\xi\right]\right)\right|_{t=0}\\[0.2em]
&=\left.\frac{d}{dt}\left[\exp(tX)e_{p+1}:\exp(tY)(w+t\xi)\right]\,\right|_{t=0}\\[0.2em]
&=(d\pi)_{(e_{p+1},w)}\left(Xe_{p+1},Yw+\xi\right).
\end{align*}

Suppose that this vector is zero. Since $\ker(d\pi)_{(e_{p+1},w)}=\mathbb{C}(e_{p+1},w),$ there exists $\mu\in\mathbb{C}$ such that $Xe_{p+1}=\mu e_{p+1},$ $Yw+\xi=\mu w.$ Since $X\in\mathfrak{u}(p+1)$, we have $\mu=ia$ for some $a\in\mathbb{R}$ and $X=\mathrm{diag}(X_1,ia),$ $X_1\in\mathfrak{u}(p).$ Hence, $A:=\left(\mathrm{diag}(X_1,ia),Y\right)\in\mathfrak{k}_0.$ Moreover,
\[
(d\sigma_0)_{(\mathrm{Id}_{p+1},\mathrm{Id}_{q+1})}(A)=Y-ia\mathrm{Id}_{q+1},
\]
and therefore
\[
\xi=iaw-Yw=-(d\sigma_0)_{(\mathrm{Id}_{p+1},\mathrm{Id}_{q+1})}(A)w.
\]
Thus $[(X,Y),\xi]=[(0,0),0],$ and consequently $\ker(d\Psi_0)=\{0\}.$ Since
\[
\dim_{\mathbb{R}}\left(G_{p,q}\times_{K_0}\mathbb{C}^{q+1}\right)=\dim_{\mathbb{R}}(G_{p,q}/K_0)+2(q+1)=2(p+q+1)=\dim_{\mathbb{R}}\mathcal{U}_0,
\]
the differential $d\Psi_0$ is an isomorphism at every point. Hence $\Psi_0$ is a local diffeomorphism. Since $\Psi_0$ is bijective, it is a diffeomorphism.\\

Now, for $r\geq0$ sufficiently small, we have
\begin{align*}
\Psi_0\left([(\mathrm{Id}_{p+1},\mathrm{Id}_{q+1}),rf_{q+1}]\right)&=[e_{p+1}:rf_{q+1}]\\[0.2em]
&=\left[\frac{1}{\sqrt{1+r^2}}e_{p+1}:\frac{r}{\sqrt{1+r^2}}f_{q+1}\right]\\[0.2em]
&=[(\cos(\arctan r))e_{p+1}:(\sin(\arctan r))f_{q+1}]\\[0.2em]
&=\gamma(\arctan r).
\end{align*}
Recall that $\tau_0(s)=s\beta_0(s^2).$ Since $\tan x$ is an odd and smooth function near 0, there exists a smooth function $\alpha$ defined near $0$ such that $\alpha(0)=1$ and $\tan x=x\alpha(x^2).$ Therefore,
\begin{align*}
\tan(\tau_0(s))&=\tau_0(s)\alpha(\tau_0(s)^2)\\
&=s\beta_0(s^2)\alpha\left(s^2\beta_0(s^2)^2\right).
\end{align*}
Defining $\psi_0(r):=\beta_0(r)\alpha\left(r\beta_0(r)^2\right),$ we obtain 
\[
\tan(\tau_0(s))=s\psi_0(s^2),
\]
where the function $\psi_0$ is smooth and
\[
\psi_0(0)=\beta_0(0)\alpha(0)=\frac{1}{f_0(0)}>0.
\]
Decreasing $\varepsilon_0$ if necessary, we may assume that $\psi_0$ is positive on $[0,\varepsilon_0^2)$. Let us define 
\[
R_0:\mathbb{C}^{q+1}\supseteq B(0;\varepsilon_0)\ni w\longmapsto \psi_0(||w||^2)w\in \mathbb{C}^{q+1}.
\]
Since the representation $\sigma_0$ is unitary, $||\sigma_0(k)w||=||w||,$ for $k\in K_0$ and $w\in\mathbb{C}^{q+1}.$ Therefore, $R_0$ is $K_0$-equivariant. Furthermore, 
\[
(dR_0)_0=\psi_0(0)\mathrm{Id}_{\mathbb{C}^{q+1}}=\frac{1}{f_0(0)}\mathrm{Id}_{\mathbb{C}^{q+1}},
\]
so, restricted to a sufficiently small neighborhood $B$ of $0$, $R_0$  is a smooth $K_0$-equivariant diffeomorphism onto its image. Consequently, $R_0$ induces a $G_{p,q}$-equivariant smooth diffeomorphism 
\[
\widehat{R}_0:G_{p,q}\times_{K_0} B\ni[a,w] \longmapsto [a,R_0(w)]\in G_{p,q}\times_{K_0}R_0(B).
\]
Define $\Phi_0:=\Psi_0\circ \widehat{R}_0.$ Then, $\Phi_0$ is a $G_{p,q}$-equivariant smooth diffeomorphism from $G_{p,q}\times_{K_0}B$ onto a neighborhood of $Q_0$ in $\mathbb{C}P^{p+q+1}$ and, for $s\geq 0$ sufficiently small,
\begin{align*}
    \Phi_0\left(\left[(\mathrm{Id}_{p+1},\mathrm{Id}_{q+1}),sf_{q+1}\right]\right)&=\Psi_0\left(\widehat{R}_0\left(\left[(\mathrm{Id}_{p+1},\mathrm{Id}_{q+1}),sf_{q+1}\right]\right)\right)\\[0.2em]
    &=\Psi_0\left(\left[(\mathrm{Id}_{p+1},\mathrm{Id}_{q+1}),s\psi_0(s^2)f_{q+1}\right]\right)\\[0.2em]
    &=\gamma\left(\arctan\left(s\psi_0(s^2)\right)\right)\\[0.2em]
    &=\gamma\left(\arctan\left(\tan(\tau_0(s)\right)\right)\\[0.2em]
    &=\gamma(\tau_0(s))\\[0.2em]
    &=c_0(s).
\end{align*}
This proves $b).$
\end{proof}
Let $\widetilde{g}_0:=\Phi_0^*g$ and set $\widetilde{c}_0(s):=\left[(\mathrm{Id}_{p+1},\mathrm{Id}_{q+1}),sf_{q+1}\right].$ By Lemma~\ref{lem:arc-length-Q0}, $\Phi_0(\widetilde{c}_0(s))=c_0(s)$ for all $s\geq 0$ sufficiently small. Since $\Phi_0$ is $G_{p,q}$-equivariant, the metric $\widetilde{g}_0$ is $G_{p,q}$-invariant and, for $X,Y\in\mathfrak{m},$ we have 
\[
\left(\widetilde{g}_0\right)_{\widetilde{c}_0(s)}(X^\#(\widetilde{c}_0(s)),Y^\#(\widetilde{c}_0(s)))=g_{c_0(s)}(X^*(c_0(s)),Y^*(c_0(s))).
\]
Furthermore, for $s>0$ sufficiently small, $\widetilde{c}_0$ is a unit-speed geodesic with respect to $\widetilde{g}_0$ and is orthogonal to every principal orbit it meets. In this context, $g$ extends to a smooth Riemannian metric on a neighborhood of $Q_0$ if and only if $\widetilde{g}_0$ extends to a smooth Riemannian metric on a neighborhood of the  zero section in $G_{p,q}\times_{K_0}\mathbb{C}^{q+1}.$
\begin{proposition}\label{prop:smooth-extension}
Let $g=(f_0,\sin(2t)f_0,f_2,f_3,0)$ and assume that $f_0$ satisfies \eqref{eq:f0-smoothness}. Then $g$ extends to a smooth Riemannian metric on a neighborhood of $Q_0$ if and only if there exist $\delta_0>0$ and smooth positive functions
\[
\phi_{2,0},\phi_{3,0}:[0,\delta_0^2)\longrightarrow\mathbb{R}
\]
such that
\begin{equation}\label{eq:smoothness-at-Q0}
\begin{cases}
&f_2(t)=\phi_{2,0}(t^2),\\[0.2em]
&f_3(t)=t\phi_{3,0}(t^2),\\[0.2em]
&\phi_{3,0}(0)=f_0(0),
\end{cases}\,\,\, 0<t<\delta_0.
\end{equation}
Similarly, $g$ extends to a smooth Riemannian metric on a neighborhood of $Q_1$ if and only if there exist $\delta_1>0$ and smooth positive functions
\[
\phi_{2,1},\phi_{3,1}:[0,\delta_1^2)\longrightarrow\mathbb{R}
\]
such that
\begin{equation}\label{eq:smoothness-at-Q1}
\begin{cases}
&f_2\left(\frac{\pi}{2}-t\right)=t\phi_{2,1}(t^2),\\[0.2em]
&f_3\left(\frac{\pi}{2}-t\right)=\phi_{3,1}(t^2),\\[0.2em]
&\phi_{2,1}(0)=f_0\left(\frac{\pi}{2}\right),
\end{cases}\,\,\, 0<t<\delta_1.
\end{equation}
\end{proposition}
\begin{proof}
First, note that $f_{q+1}\in\mathbb{S}^{2q+1}\subseteq\mathbb{C}^{q+1}$ and that the isotropy group of $f_{q+1}$ in $K_0$ with respect to the representation $\sigma_0$ is $H.$ Therefore, we can apply Proposition~\ref{prop:V-Z-smoothness-conditions} to $\widetilde{g}_0.$ Recall that $\mathfrak{p}_0=\mathfrak{m}_1\oplus\mathfrak{m}_3,$ $\mathfrak{n}_0=\mathfrak{m}_2.$ Let us choose the vectors $Z\in\mathfrak{m}_1$ and $C_q\in\mathfrak{m}_3.$ For each $\theta \in\mathbb{R},$ we have
\begin{align*}
    &\exp(\theta Z)=\left(\mathrm{diag}\left(\mathrm{Id}_p,e^{i\theta}\right),\mathrm{diag}\left(\mathrm{Id}_q,e^{-i\theta}\right)\right),\\[0.2em]
    &\exp(\theta C_q)=\left(\mathrm{Id}_{p+1},\mathrm{diag}\left(\mathrm{Id}_{q-1},\begin{pmatrix}\cos\theta&-\sin\theta\\\sin\theta&\cos\theta\end{pmatrix}\right)\right).
\end{align*}
Therefore, 
\[
L_Z:=\{\exp\left(\theta Z\right):\theta\in\mathbb{R}\},\ L_{C_q}:=\{\exp(\theta C_q):\theta\in\mathbb{R}\}
\]
are closed one-parameter subgroups of $K_0$ and both maps $\theta \longmapsto\exp(\theta Z),$ $\theta\longmapsto \exp(\theta C_q)$ have period $2\pi.$ Let us compute $|L_Z\cap H|$ and $|L_{C_q}\cap H|.$ An element $\exp(\theta Z)$ belongs to $H$ if and only if $e^{i\theta}=e^{-i\theta}.$ Thus, $e^{2i\theta}=1$ and hence $\theta\in\pi\mathbb{Z}$. Since the map $\theta\mapsto\exp(\theta Z)$ has period $2\pi$, we obtain
\[
a_Z:=|L_Z\cap H|=2.
\]
On the other hand, if $\exp(\theta C_q)\in H$, then there exist $U_1\in\mathrm{U}(p),$ $\lambda\in\mathrm{U}(1)$ and $V_1\in\mathrm{U}(q)$ such that
\begin{align*}
\mathrm{Id}_{p+1}=\mathrm{diag}(U_1,\lambda),\ \mathrm{diag}(\mathrm{Id}_{q-1},\cos\theta)=V_1,\ \cos\theta=\lambda,\ \sin\theta=0.
\end{align*}
Therefore, $U_1=\mathrm{Id}_p,$ $\lambda=1,$ $V_1=\mathrm{Id}_{q}$ and $\theta\in2\pi\mathbb{Z},$ that is $L_{C_q}\cap H=\{(\mathrm{Id}_{p+1},\mathrm{Id}_{q+1})\}.$ Consequently,
\[
a_{C_q}:=|L_{C_q}\cap H|=1.
\]
For $A_j,B_j\in\mathfrak{m}_2,\ j=1,...,p,$ we have
\[
[Z,A_j]=B_j,\ [Z,B_j]=-A_j,\ j=1,...,p,
\]
and $L_{C_q}$ acts trivially on $\mathfrak{m}_2$ because $[C_q,A_j]=[C_q,B_j]=0,\ j=1,...,p.$ By Proposition~\ref{prop:V-Z-smoothness-conditions} and the discussion after the proof of Lemma~\ref{lem:arc-length-Q0}, we have that $g$ extends to a smooth $G_{p,q}$-invariant symmetric $2$-tensor on a neighborhood of $Q_0$ if and only if
\begin{align*}
&g_{c_0(s)}\left(Z^*(c_0(s)),Z^*(c_0(s))\right)=4s^2+s^4\phi_Z(s^2),\\
&g_{c_0(s)}\left(A_j^*(c_0(s)),A_j^*(c_0(s))\right)+g_{c_0(s)}\left(B_j^*(c_0(s)),B_j^*(c_0(s))\right)=\phi_{1,Z,j}(s^2),\\
&g_{c_0(s)}\left(A_j^*(c_0(s)),A_j^*(c_0(s))\right)-g_{c_0(s)}\left(B_j^*(c_0(s)),B_j^*(c_0(s))\right)=s\phi_{2,Z,j}(s^2),\\
&g_{c_0(s)}\left(A_j^*(c_0(s)),B_j^*(c_0(s))\right)=s\phi_{3,Z,j}(s^2),\ j=1,...,p,\\
&g_{c_0(s)}\left(A_j^*(c_0(s)),A_k^*(c_0(s))\right)+g_{c_0(s)}\left(B_j^*(c_0(s)),B_k^*(c_0(s))\right)=\phi_{1,Z,j,k}(s^2),\\
&g_{c_0(s)}\left(A_j^*(c_0(s)),B_k^*(c_0(s))\right)-g_{c_0(s)}\left(B_j^*(c_0(s)),A_k^*(c_0(s))\right)=\phi_{2,Z,j,k}(s^2),\\
&g_{c_0(s)}\left(A_j^*(c_0(s)),A_k^*(c_0(s))\right)-g_{c_0(s)}\left(B_j^*(c_0(s)),B_k^*(c_0(s))\right)=s\phi_{3,Z,j,k}(s^2),\\
&g_{c_0(s)}\left(A_j^*(c_0(s)),B_k^*(c_0(s))\right)+g_{c_0(s)}\left(B_j^*(c_0(s)),A_k^*(c_0(s))\right)=s\phi_{4,Z,j,k}(s^2),\ 1\leq j<k\leq p,\\
&g_{c_0(s)}\left(C_q^*(c_0(s)),C_q^*(c_0(s))\right)=s^2+s^4\phi_{C_q}(s^2),\\
&g_{c_0(s)}\left(A_j^*(c_0(s)),A_k^*(c_0(s))\right)=\phi^{A}_{j,k}(s^2),\ 1\leq j\leq k\leq p,\\
&g_{c_0(s)}\left(B_j^*(c_0(s)),B_k^*(c_0(s))\right)=\phi^{B}_{j,k}(s^2),\ 1\leq j\leq k\leq p,\\
&g_{c_0(s)}\left(A_j^*(c_0(s)),B_k^*(c_0(s))\right)=\phi^{AB}_{j,k}(s^2),\ 1\leq j,k\leq p,
\end{align*}
for some functions $\phi_Z,\phi_{\ell,Z,j},\phi_{\ell,Z,j,k},\phi_{C_q},\phi^A_{j,k},\phi^B_{j,k},\phi^{AB}_{j,k}$ that are smooth on an interval of the form $[0,\tilde{\delta}^2),$ for some $\tilde{\delta}>0.$ Since $c_0(s)=\gamma(\tau_0(s))$ and by the definition of $g$, all the conditions above reduce to
\[
\begin{cases}
&\sin^2\left(2\tau_0(s)\right)f_0(\tau_0(s))^2
=4s^2+s^4\psi_1(s^2),\\
&f_2(\tau_0(s))^2=\psi_2(s^2),\\
&f_3(\tau_0(s))^2=s^2+s^4\psi_3(s^2),
\end{cases}\ 0<s<\widetilde{\delta}
\]
where $\psi_1,\psi_2,\psi_3$ are smooth functions on an interval of the form $[0,\tilde{\delta}^2).$ The resulting extension is a smooth $G_{p,q}$-invariant Riemannian metric if and only if $\psi_2(0)>0$, since $\mathfrak{n}_0=\mathfrak{m}_2$ is tangent to $Q_0$ and the restriction of the extension to $\mathfrak{m}_2$ is $\psi_2(0)(\cdot,\cdot)$. We will show that the first condition is a consequence of \eqref{eq:f0-smoothness}. In fact, by the proof of Lemma~\ref{lem:arc-length-Q0}, we have that $\tau_0(s)=s\beta_0(s^2)$ for some smooth function $\beta_0$ with $\beta_0(0)=1/f_0(0).$ Since $\sin(2x)$ is smooth and odd, $\sin(2x)=2x\xi(x^2),$ where $\xi$ is smooth and $\xi(0)=1.$ Therefore,
\begin{align*}
    \sin^2(2\tau_0(s))f_0(\tau_0(s))^2&=\left(2\tau_0(s)\xi(\tau_0(s)^2)\right)^2\phi_0\left(\tau_0(s)^2\right)^2\\[0.2em]
    &=4s^2\beta_0(s^2)^2\xi\left(s^2\beta_0(s^2)^2\right)^2\phi_0\left(s^2\beta_0(s^2)^2\right)^2.
\end{align*}
Defining $\kappa(r):=\beta_0(r)^2\xi\left(r\beta_0(r)^2\right)^2\phi_0\left(r\beta_0(r)^2\right)^2,$ we have that $\kappa$ is smooth and 
\[
\kappa(0)=\beta_0(0)^2\xi(0)^2\phi_0(0)^2=(1/f_0(0))^2\cdot(1)\cdot (f_0(0))^2=1.
\]
Hence,
\[
\displaystyle\kappa(r)=1+r\int_0^1\kappa'(ur)\,du=1+r\chi(r),
\]
where $\displaystyle\chi(r):=\int_0^1\kappa'(ur)\,du$ is smooth. Therefore,
\begin{align*}
    \sin^2(2\tau_0(s))f_0(\tau_0(s))^2&=4s^2\kappa(s^2)=4s^2+4s^4\chi(s^2)=4s^2+s^4\psi_1(s^2),
\end{align*}
with $\psi_1(r):=4\chi(r)$ smooth. Thus, under \eqref{eq:f0-smoothness}, the first condition imposes no further restrictions. Consequently, $g$ extends to a smooth Riemannian metric on a neighborhood of $Q_0$ if and only if there exist smooth functions $\psi_2,\psi_3$ on an interval of the form $[0,\widetilde{\delta}^2)$ such that
\begin{equation}\label{eq:extension-conditions-s}
\begin{cases}
    &f_2(\tau_0(s))^2=\psi_2(s^2),\\
    &f_3(\tau_0(s))^2=s^2+s^4\psi_3(s^2),\\
    &\psi_2(0)>0,
\end{cases}\ 0<s<\widetilde{\delta}.
\end{equation}
We will show that \eqref{eq:extension-conditions-s} is equivalent to \eqref{eq:smoothness-at-Q0}. Note that
\[
s_0(t)=\int_0^t f_0(r)\,dr=\int_0^t\phi_0(r^2)\,dr=t\int_0^1\phi_0(u^2t^2)\,du=t\mu_0(t^2),
\]
where $\displaystyle\mu_0(r):=\int_0^1\phi_0(u^2r)\,du$ is smooth and positive, and $\mu_0(0)=\phi_0(0)=f_0(0).$ Assume that \eqref{eq:extension-conditions-s} holds. Since $\psi_2(0)>0,$ $\sqrt{\psi_2}$ is smooth and positive, and since $s_0,\tau_0$ are inverse to each other,
\begin{align*}
    f_2(t)=f_2(\tau_0(s_0(t)))=\sqrt{\psi_2(s_0(t)^2)}=\sqrt{\psi_2\left(t^2\mu_0(t^2)^2\right)}.
\end{align*}
Therefore, $f_2(t)=\phi_{2,0}(t^2),$ where $\phi_{2,0}(r):=\sqrt{\psi_2\left(r\mu_0(r)^2\right)}$ is smooth and positive. On the other hand,
\[
f_3(\tau_0(s))^2=s^2+s^4\psi_3(s^2)=s^2(1+s^2\psi_3(s^2)).
\]
Since $1+r\psi_3(r)$ is smooth and positive, $\sqrt{1+r\psi_3(r)}$ is smooth and positive. It follows that
\begin{align*}
    f_3(t)&=f_3(\tau_0(s_0(t)))\\[0.2em]
          &=s_0(t)\sqrt{1+s_0(t)^2\psi_3(s_0(t)^2)}\\[0.2em]
          &=t\mu_0(t^2)\sqrt{1+t^2\mu_0(t^2)^2\psi_3(t^2\mu_0(t^2)^2)}\\[0.2em]
          &=t\phi_{3,0}(t^2),
\end{align*}
where $\phi_{3,0}(r):=\mu_0(r)\sqrt{1+r\mu_0(r)^2\psi_3\left(r\mu_0(r)^2\right)}$ is smooth and positive, and $\phi_{3,0}(0)=\mu_0(0)=f_0(0).$ Therefore, we have \eqref{eq:smoothness-at-Q0}.\\

Conversely, suppose that \eqref{eq:smoothness-at-Q0} holds. Since $\tau_0(s)=s\beta_0(s^2),$ we obtain
\[
f_2(\tau_0(s))^2=\phi_{2,0}\left(s^2\beta_0(s^2)^2\right)^2=\psi_2(s^2),
\]
where $\psi_2(r):=\phi_{2,0}\left(r\beta_0(r)^2\right)^2$ is smooth and $\psi_2(0)=\phi_{2,0}(0)^2>0.$ Also
\begin{align*}
f_3(\tau_0(s))=s\beta_0(s^2)\phi_{3,0}\left(s^2\beta_0(s^2)^2\right)=s\eta(s^2),
\end{align*}
where $\eta(r):=\beta_0(r)\phi_{3,0}\left(r\beta_0(r)^2\right).$ The function $\eta$ is smooth and $\eta(0)=\beta_0(0)\phi_{3,0}(0)=\frac{1}{f_0(0)}f_0(0)=1.$ Hence, there exists a smooth function $\psi_3$ such that $\eta(r)^2=1+r\psi_3(r).$ Therefore,
\[
f_3(\tau_0(s))^2=s^2\eta(s^2)^2=s^2+s^4\psi_3(s^2).
\]
Thus, \eqref{eq:extension-conditions-s} holds.\\

The statement for $Q_1$ follows analogously. Set $g_j(t):=f_j\left(\frac{\pi}{2}-t\right),\ j=0,2,3.$ Using the analogue of Lemma~\ref{lem:arc-length-Q0} at $Q_1$ and $\mathfrak{p}_1=\mathfrak{m}_1\oplus\mathfrak{m}_2,$ $\mathfrak{n}_1=\mathfrak{m}_3,$ the same argument applied to $g_0,g_2,g_3$ gives the conditions at $Q_1$.
\end{proof}
\begin{example}[The Fubini-Study metric]\label{ex:Fubini-Study}
Let $g_{\mathrm{FS}}$ be the Fubini-Study metric on $\mathbb{C}P^{p+q+1}$ induced by the standard metric on the unit sphere $\mathbb{S}^{2p+2q+3}\subseteq\mathbb{C}^{p+1}\times\mathbb{C}^{q+1}.$ We have that
\begin{equation*}\label{eq:Fubini-Study-metric}
    g_{\mathrm{FS}}=(1,\sin(2t),\cos t,\sin t,0).
\end{equation*}
To verify this, consider the canonical projection
\[
\pi\big{|}_{\mathbb{S}^{2p+2q+3}}:\mathbb{S}^{2p+2q+3}\ni(z,w)
\longrightarrow[z:w]\in\mathbb{C}P^{p+q+1}.
\]
The curve $\gamma$ is the projection of $z(t)=\left((\cos t)e_{p+1},(\sin t)f_{q+1}\right).$ Since $z'(t)=\left(-(\sin t)e_{p+1},(\cos t)f_{q+1}\right)$ has unit length and is orthogonal to $iz(t)$, we have $g_{\mathrm{FS}}(\gamma'(t),\gamma'(t))=1,$ so $f_0(t)=1.$ For $j=1,...,p$, the vectors induced by $A_j$ and $B_j$ at $z(t)$ are $(-\cos(t)e_j,0)$ and $(i\cos(t)e_j,0),$ respectively. They are orthogonal to $iz(t),$ therefore
\[
\left(g_{\mathrm{FS}}\right)_{\gamma(t)}(A_j^*(\gamma(t)),A_j^*(\gamma(t)))=\left(g_{\mathrm{FS}}\right)_{\gamma(t)}(B_j^*(\gamma(t)),B_j^*(\gamma(t)))=\cos^2t.
\]
Thus, $f_2(t)=\cos t.$ Similarly, for $j=1,...,q$, the vectors induced by $C_j$ and $D_j$ are $(0,-\sin(t)f_j)$ and $(0,i\sin(t)f_j),$ so
\[
\left(g_{\mathrm{FS}}\right)_{\gamma(t)}(C_j^*(\gamma(t)),C_j^*(\gamma(t)))=\left(g_{\mathrm{FS}}\right)_{\gamma(t)}(D_j^*(\gamma(t)),D_j^*(\gamma(t)))=\sin^2t.
\]
Therefore, $f_3(t)=\sin t.$ Since the Fubini-Study metric is compatible with the standard complex structure $J$, Corollary~\ref{coro:compatible-metrics} implies that $u=0,$ $f_1(t)=\sin(2t)f_0(t)=\sin(2t).$ Therefore, $g_{\mathrm{FS}}=(1,\sin(2t),\cos t,\sin t,0).$
\end{example}

\section{Balanced $G_{p,q}$-invariant metrics}\label{section:balanced-metrics}
In this section, we consider a $G_{p,q}$-invariant Riemannian metric $g=(f_0,\sin(2t)f_0,f_2,f_3,0)$ on $\mathbb{C}P^{p+q+1},$ compatible with $J$. By the results of the previous section, there exist smooth positive functions $\phi_k,\phi_{j,k},\ j=2,3,\ k=0,1,$ defined near $0,$ satisfying the equations \eqref{eq:f0-smoothness}, \eqref{eq:smoothness-at-Q0} and \eqref{eq:smoothness-at-Q1}. Let $t:\left(\mathbb{C}P^{p+q+1}\right)_0\longrightarrow\left(0,\frac{\pi}{2}\right)$ be the $G_{p,q}$-invariant function defined by
\[
t([z:w])=t_0\Longleftrightarrow[z:w]\in G_{p,q}\cdot\gamma(t_0).
\]
As observed in Subsection~\ref{subsec:invariant-metrics}, $dt$ is a $G_{p,q}$-invariant $1$-form on $\left(\mathbb{C}P^{p+q+1}\right)_0$ satisfying
\[
(dt)_{\gamma(t)}(\gamma'(t))=1,\ (dt)_{\gamma(t)}(X^*(\gamma(t)))=0,\ X\in\mathfrak{m}.
\]
Let $\zeta$ be the $G_{p,q}$-invariant $1$-form on $\left(\mathbb{C}P^{p+q+1}\right)_0$ determined along $\gamma$ by
\[
\zeta_{\gamma(t)}(Z^*(\gamma(t)))=1,\ \zeta_{\gamma(t)}(\gamma'(t))=0,\ \zeta_{\gamma(t)}\left(X^*(\gamma(t))\right)=0,\ X\in\mathfrak{m}_2\oplus\mathfrak{m}_3.
\]
We also define the $G_{p,q}$-invariant $2$-forms $\omega_2$ and $\omega_3$ on $\left(\mathbb{C}P^{p+q+1}\right)_0$ by declaring, along $\gamma$,
\begin{align*}
&\left(\omega_2\right)_{\gamma(t)}(X^*(\gamma(t)),Y^*(\gamma(t)))=(JX^*(\gamma(t)),Y^*(\gamma(t)))_t,\ \textnormal{if}\ X,Y\in\mathfrak{m}_2,\\
&\left(\omega_2\right)_{\gamma(t)}(u,v)=0,\ \text{if}\ u\notin \mathfrak{m}_2\ \text{or}\ v\notin \mathfrak{m}_2,\\
&\left(\omega_3\right)_{\gamma(t)}(X^*(\gamma(t)),Y^*(\gamma(t)))=(JX^*(\gamma(t)),Y^*(\gamma(t)))_t,\ X,Y\in\mathfrak{m}_3,\\
&\left(\omega_3\right)_{\gamma(t)}(u,v)=0,\ \text{if}\ u\notin \mathfrak{m}_3\ \text{or}\ v\notin \mathfrak{m}_3.
\end{align*}
Recall that the fundamental form $\omega$ of the Hermitian structure $(g,J)$ is defined by
\[
\omega(X,Y):=g(JX,Y).
\]
We first describe $\omega$ on $\left(\mathbb{C}P^{p+q+1}\right)_0$.
\begin{proposition}\label{prop:fundamental-form}
The fundamental form of $(g,J)$ on $\left(\mathbb{C}P^{p+q+1}\right)_0$ is
\[
\omega=-\sin(2t)f_0(t)^2\,dt\wedge\zeta+f_2(t)^2\omega_2+f_3(t)^2\omega_3.
\]
\end{proposition}
\begin{proof}
Since $\omega$, $dt$, $\zeta$, $\omega_2$ and $\omega_3$ are $G_{p,q}$-invariant, it is enough to compute $\omega$ at $\gamma(t),\ 0<t<\pi/2.$ By Proposition~\ref{prop:canonical-J},
\[
J\gamma'(t)=-\frac{1}{\sin(2t)}Z^*(\gamma(t)).
\]
Therefore,
\begin{align*}
\omega_{\gamma(t)}\left(\gamma'(t),Z^*(\gamma(t))\right)&=g_{\gamma(t)}\left(J\gamma'(t),Z^*(\gamma(t))\right)\\[0.2em]
&=-\frac{1}{\sin(2t)}g_{\gamma(t)}\left(Z^*(\gamma(t)),Z^*(\gamma(t))\right)\\
&=-\sin(2t)f_0(t)^2.
\end{align*}

But $(dt\wedge\zeta)_{\gamma(t)}\left(\gamma'(t),Z^*(\gamma(t))\right)=1.$ Hence, the restriction of $\omega$ to $\mathbb{R}\gamma'(t)\oplus\mathfrak{m}_1$ is 
\[
-\sin(2t)f_0(t)^2\,dt\wedge\zeta.
\]
For $X,Y\in\mathfrak{m}_2$, we have
\begin{align*}
\omega_{\gamma(t)}(X^*(\gamma(t)),Y^*(\gamma(t)))&=g_{\gamma(t)}(JX^*(\gamma(t)),Y^*(\gamma(t)))\\[0.2em]
&=f_2(t)^2(JX^*(\gamma(t)),Y^*(\gamma(t)))_t\\[0.2em]
&=f_2(t)^2(\omega_2)_{\gamma(t)}(X^*(\gamma(t)),Y^*(\gamma(t))).
\end{align*}
Similarly, for $X,Y\in\mathfrak{m}_3$,
\[
\omega_{\gamma(t)}(X^*(\gamma(t)),Y^*(\gamma(t)))=f_3(t)^2(\omega_3)_{\gamma(t)}(X^*(\gamma(t)),Y^*(\gamma(t))).
\]
The result then follows from the fact that $J$ preserves $\mathbb{R}\gamma'(t)\oplus\mathfrak{m}_1,$ $\mathfrak{m}_2,$ $\mathfrak{m}_3,$ and these subspaces are pairwise orthogonal with respect to $g.$
\end{proof}
Recall that a Hermitian structure $(g,J)$ on a complex manifold of complex dimension $n$ is called K\"ahler if its fundamental form $\omega$ satisfies $d\omega=0$, and balanced if $d\omega^{n-1}=0$ (see \cite[Theorem~1.6]{M1982}). In our case, $n=p+q+1$. Since $g$ and $J$ are smooth on $\mathbb{C}P^{p+q+1}$, the forms $d\omega$ and $d\omega^{p+q}$ are smooth on $\mathbb{C}P^{p+q+1}$. Therefore, as $\left(\mathbb{C}P^{p+q+1}\right)_0$ is dense, $(g,J)$ is K\"ahler if and only if $d\omega=0$ on $\left(\mathbb{C}P^{p+q+1}\right)_0$, and balanced if and only if $d\omega^{p+q}=0$ on $\left(\mathbb{C}P^{p+q+1}\right)_0$.
\begin{lemma}\label{lem:differential-forms}
On $\left(\mathbb{C}P^{p+q+1}\right)_0$, we have
\begin{itemize}
    \item[$a)$] $d\zeta=\omega_2-\omega_3.$
    
    \item[$b)$] $d\omega_2=d\omega_3=0.$
\end{itemize}
\end{lemma}
\begin{proof} Since $\zeta$, $\omega_2$ and $\omega_3$ are $G_{p,q}$-invariant, it is enough to verify the identities at $\gamma(t)$, $0<t<\pi/2$. Consider the diffeomorphism
\[
\Psi:\left(0,\frac{\pi}{2}\right)\times G_{p,q}/H\ni(t,(U,V)H)\longmapsto (U,V)\cdot\gamma(t)\in\left(\mathbb{C}P^{p+q+1}\right)_0
\]
and set 
\begin{align*}
    \widetilde{\gamma'}:=\Psi_*\left(\partial_t\right),\ \widehat{X}(t,(U,V)H):=\left.\left(0,\frac{d}{ds}\exp\left(sX\right)(U,V)H\,\right|_{s=0}\right),\ X\in\mathfrak{m}.
\end{align*}
Then, $\widetilde{\gamma'}(\gamma(t))=\gamma'(t),\ \Psi_*\widehat{X}=X^*,$ and $[\widetilde{\gamma'},X^*]=0,$ for all $X\in\mathfrak{m}.$ Note that $\widetilde{\gamma'}$ is a $G_{p,q}$-invariant vector field on $(\mathbb{C}P^{p+q+1})_0.$ Given $X\in\mathfrak{m},$ we have
\begin{align*}
    (d\zeta)_{\gamma(t)}(\gamma'(t),X^*(\gamma(t)))&=\gamma'(t)(\zeta(X^*))-X^*(\gamma(t))(\zeta(\widetilde{\gamma'}))-\zeta_{\gamma(t)}\left([\widetilde{\gamma'},X^*](\gamma(t))\right)\\[0.2em]
    &=\gamma'(t)(\zeta(X^*))-X^*(\gamma(t))(\zeta(\widetilde{\gamma'})).
\end{align*}
By definition, $\zeta(X^*)(\gamma(t))=\zeta_{\gamma(t)}(X^*(\gamma(t)))$ is independent of $t,$ and, since $\zeta$ and $\widetilde{\gamma'}$ are $G_{p,q}$-invariant, $\zeta(\widetilde{\gamma'})(\exp(sX)\cdot\gamma(t))=\zeta_{\gamma(t)}(\gamma'(t))=0.$ Thus,
\[
(d\zeta)_{\gamma(t)}(\gamma'(t),X^*(\gamma(t)))=0.
\]
Also, for $X,Y\in\mathfrak{m}$ we have
\begin{align*}
    (d\zeta)_{\gamma(t)}(X^*(\gamma(t)),Y^*(\gamma(t)))&=X^*(\gamma(t))(\zeta(Y^*))-Y^*(\gamma(t))(\zeta(X^*))-\zeta_{\gamma(t)}\left([X^*,Y^*](\gamma(t))\right)\\[0.2em]
    &=X^*(\gamma(t))(\zeta(Y^*))-Y^*(\gamma(t))(\zeta(X^*))+\zeta_{\gamma(t)}\left([X,Y]_{\mathfrak{m}}^*(\gamma(t))\right).
\end{align*}
But 
\begin{align*}
\zeta(Y^*)(\exp(sX)\cdot \gamma(t))&=\zeta_{\exp(sX)\cdot \gamma(t)}\left(Y^*\left(\exp(sX)\cdot \gamma(t)\right)\right)\\[0.2em]
&=\zeta_{\exp(sX)\cdot \gamma(t)}\left(\left(d\phi_{\exp(sX)}\right)_{\gamma(t)}\left(\mathrm{Ad}(\exp(-sX))Y\right)^*(\gamma(t))\right)\\[0.2em]
&=\zeta_{\gamma(t)}\left(\left(\mathrm{Ad}(\exp(-sX))Y\right)^*(\gamma(t))\right),
\end{align*}
and
\begin{align*}
X^*(\gamma(t))(\zeta(Y^*))&=\left.\frac{d}{ds}\zeta(Y^*)(\exp(sX)\cdot\gamma(t))\,\right|_{s=0}\\[0.5em]
&=\left.\frac{d}{ds}\zeta_{\gamma(t)}\left(\left(\mathrm{Ad}(\exp(-sX))Y\right)^*(\gamma(t))\right)\,\right|_{s=0}\\[0.5em]
&=\zeta_{\gamma(t)}\left.\left(\frac{d}{ds}\left(\mathrm{Ad}(\exp(-sX))Y\right)^*(\gamma(t))\,\right|_{s=0}\right)\\[0.5em]
&=-\zeta_{\gamma(t)}\left([X,Y]_{\mathfrak{m}}^*(\gamma(t))\right).
\end{align*}
Hence,
\begin{align*}
(d\zeta)_{\gamma(t)}\left(X^*(\gamma(t)),Y^*(\gamma(t))\right)&=-\zeta_{\gamma(t)}\left([X,Y]_{\mathfrak{m}}^*(\gamma(t))\right)+\zeta_{\gamma(t)}\left([Y,X]_{\mathfrak{m}}^*(\gamma(t))\right)+\zeta_{\gamma(t)}\left([X,Y]_{\mathfrak{m}}^*(\gamma(t))\right)\\[0.2em]
&=-\zeta_{\gamma(t)}\left([X,Y]_{\mathfrak{m}}^*(\gamma(t))\right).
\end{align*}
It is a straightforward computation to verify that 
\[
[\mathfrak{m}_1,\mathfrak{m}_2]\subseteq\mathfrak{m}_2,\ [\mathfrak{m}_1,\mathfrak{m}_3]\subseteq\mathfrak{m}_3,\ [\mathfrak{m}_2,\mathfrak{m}_3]=\{0\},\\[0.2em]
[\mathfrak{m}_2,\mathfrak{m}_2]\subseteq\mathfrak{h}\oplus\mathfrak{m}_1,\ [\mathfrak{m}_3,\mathfrak{m}_3]\subseteq\mathfrak{h}\oplus\mathfrak{m}_1
\]
and that the only brackets that may be nonzero among the vectors $Z,A_j,B_j,C_k,D_k,\ j=1,...p,\ k=1,...,q,$ are
\begin{align*}
    &[Z,A_j]=B_j,\ [Z,B_j]=-A_j,\ j=1,...,p,\\
    &[Z,C_j]=-D_j,\ [Z,D_j]=C_j,\ j=1,...,q,\\
    &[A_j,A_k]=[B_j,B_k]=(E_{k,j}-E_{j,k},0)\in\mathfrak{h},\\
    &[C_j,C_k]=[D_j,D_k]=(0,F_{k,j}-F_{j,k})\in\mathfrak{h},\\
    &[A_j,B_k]=\delta_{j}^kZ+(i(\delta_j^kE_{p+1,p+1}-E_{j,k}-E_{k,j}),i\delta_j^kF_{q+1,q+1}),\ j,k=1,...,p,\\
    &[C_j,D_k]=-\delta_{j}^kZ+(i\delta_j^kE_{p+1,p+1},i(\delta_j^kF_{q+1,q+1}-F_{j,k}-F_{k,j})),\ j,k=1,...,q,
\end{align*}
where $\delta_j^k$ is the Kronecker delta. Therefore, up to skew-symmetry, the only nonzero values  $(d\zeta)_{\gamma(t)}(u,v),\ u,v\in\mathcal{B}_t$ are
\begin{align*}
    &(d\zeta)_{\gamma(t)}\left(A_j^*(\gamma(t)),B_j^*(\gamma(t))\right)=-\zeta_{\gamma(t)}\left([A_j,B_j]_{\mathfrak{m}}^*(\gamma(t))\right)=-\zeta_{\gamma(t)}\left(Z^*(\gamma(t))\right)=-1,\ j=1,...,p,\\
    &(d\zeta)_{\gamma(t)}\left(C_j^*(\gamma(t)),D_j^*(\gamma(t))\right)=-\zeta_{\gamma(t)}\left([C_j,D_j]_{\mathfrak{m}}^*(\gamma(t))\right)=-\zeta_{\gamma(t)}\left(-Z^*(\gamma(t))\right)=1,\ j=1,...,q.    
\end{align*}
On the other hand, by the definition of $\omega_2$ and $\omega_3,$ up to skew-symmetry, the only nonzero values $(\omega_2-\omega_3)_{\gamma(t)}(u,v),\ u,v\in\mathcal{B}_t$ are
\begin{align*}
    &\begin{aligned}
    (\omega_2-\omega_3)_{\gamma(t)}(A_j^*(\gamma(t)),B_j^*(\gamma(t)))&=(\omega_2)_{\gamma(t)}(A_j^*(\gamma(t)),B_j^*(\gamma(t)))\\
    &=(JA_j^*(\gamma(t)),B_j^*(\gamma(t)))_t\\
    &=-(B_j^*(\gamma(t)),B_j^*(\gamma(t)))_t\\
    &=-1,
    \end{aligned}\ j=1,...,p,\\[0.2em]
    &\begin{aligned}
    (\omega_2-\omega_3)_{\gamma(t)}(C_j^*(\gamma(t)),D_j^*(\gamma(t)))&=-(\omega_3)_{\gamma(t)}(C_j^*(\gamma(t)),D_j^*(\gamma(t)))\\
    &=-(JC_j^*(\gamma(t)),D_j^*(\gamma(t)))_t\\
    &=(D_j^*(\gamma(t)),D_j^*(\gamma(t)))_t\\
    &=1,
    \end{aligned}\ j=1,...,q.    
\end{align*}
Hence, $(d\zeta)_{\gamma(t)}=(\omega_2-\omega_3)_{\gamma(t)}.$ This proves $a).$\\

To prove $b),$ recall that, given $X,Y\in\mathfrak{m},$ we have
\begin{align*}
    (d\omega_2)_{\gamma(t)}(\gamma'(t),X^*(\gamma(t)),Y^*(\gamma(t)))=&\,\gamma'(t)\left(\omega_2(X^*,Y^*)\right)-X^*(\gamma(t))(\omega_2(\widetilde{\gamma'},Y^*))\\[0.2em]
    &+Y^*(\gamma(t))(\omega_2(\widetilde{\gamma'},X^*))-(\omega_2)_{\gamma(t)}\left([\widetilde{\gamma'},X^*](\gamma(t)),Y^*(\gamma(t))\right)\\[0.2em]
    &+(\omega_2)_{\gamma(t)}\left([\widetilde{\gamma'},Y^*](\gamma(t)),X^*(\gamma(t))\right)\\[0.2em]
    &-(\omega_2)_{\gamma(t)}\left([X^*,Y^*](\gamma(t)),\gamma'(t)\right).
\end{align*}
Since $[\widetilde{\gamma'},X^*]=[\widetilde{\gamma'},Y^*]=0,$ $\omega_2(\widetilde{\gamma'},\cdot)=0$ and the function $t\longmapsto (\omega_2)_{\gamma(t)}(X^*(\gamma(t)),Y^*(\gamma(t)))$ is constant, we conclude
\[
(d\omega_2)_{\gamma(t)}(\gamma'(t),X^*(\gamma(t)),Y^*(\gamma(t)))=0.
\]
Now, let $X,Y,W\in\mathfrak{m}.$ Then,
\begin{align*}
    \omega_2(Y^*,W^*)(\exp(sX)\cdot \gamma(t))&=(\omega_2)_{\exp(sX)\cdot\gamma(t)}(Y^*(\exp(sX)\cdot \gamma(t)),W^*(\exp(sX)\cdot \gamma(t)))\\[0.2em]
    &=(\omega_2)_{\gamma(t)}\left((\mathrm{Ad}(\exp(-sX))Y)^*(\gamma(t)),(\mathrm{Ad}(\exp(-sX))W)^*(\gamma(t))\right).
\end{align*}
It follows that
\begin{align*}
    X^*(\gamma(t))(\omega_2(Y^*,W^*))&=\left.\frac{d}{ds}(\omega_2)_{\gamma(t)}\left((\mathrm{Ad}(\exp(-sX))Y)^*(\gamma(t)),(\mathrm{Ad}(\exp(-sX))W)^*(\gamma(t))\right)\,\right|_{s=0}\\[0.5em]
    &=-(\omega_2)_{\gamma(t)}\left([X,Y]_{\mathfrak{m}}^*(\gamma(t)),W^*(\gamma(t))\right)-(\omega_2)_{\gamma(t)}\left(Y^*(\gamma(t)),[X,W]_{\mathfrak{m}}^*(\gamma(t))\right).
\end{align*}
Similarly,
\begin{align*}
&Y^*(\gamma(t))\left(\omega_2(X^*,W^*)\right)=-(\omega_2)_{\gamma(t)}\left([Y,X]_{\mathfrak{m}}^*(\gamma(t)),W^*(\gamma(t))\right)-(\omega_2)_{\gamma(t)}\left(X^*(\gamma(t)),
[Y,W]_{\mathfrak{m}}^*(\gamma(t))\right),\\[0.2em]
&W^*(\gamma(t))\left(\omega_2(X^*,Y^*)\right)=
-(\omega_2)_{\gamma(t)}\left([W,X]_{\mathfrak{m}}^*(\gamma(t)),Y^*(\gamma(t))\right)-(\omega_2)_{\gamma(t)}\left(X^*(\gamma(t)),
[W,Y]_{\mathfrak{m}}^*(\gamma(t))\right).
\end{align*}
Evaluating at $\gamma(t)$ the identity
\begin{align*}
(d\omega_2)(X^*,Y^*,W^*)=&\,X^*\left(\omega_2(Y^*,W^*)\right)-Y^*\left(\omega_2(X^*,W^*)\right)+W^*\left(\omega_2(X^*,Y^*)\right)\\
&+\omega_2([X,Y]^*,W^*)-\omega_2([X,W]^*,Y^*)
+\omega_2([Y,W]^*,X^*),
\end{align*}
we obtain
\begin{align*}
(d\omega_2)_{\gamma(t)}\left(X^*(\gamma(t)),Y^*(\gamma(t)),W^*(\gamma(t))\right)=&\,-(\omega_2)_{\gamma(t)}\left([X,Y]_{\mathfrak{m}}^*(\gamma(t)),W^*(\gamma(t))\right)\\
&+(\omega_2)_{\gamma(t)}\left([X,W]_{\mathfrak{m}}^*(\gamma(t)),Y^*(\gamma(t))\right)\\
&-(\omega_2)_{\gamma(t)}\left([Y,W]_{\mathfrak{m}}^*(\gamma(t)),X^*(\gamma(t))\right).
\end{align*}
By the definition of $\omega_2,$ the only possible nonzero values for this expression with $X,Y,W\in\{Z,A_j,B_j,C_k,D_k:j=1,...p,\ k=1,...,q\}$ are obtained whenever $[X,Y]_{\mathfrak{m}},W\in\mathfrak{m}_2,$ $[X,W]_{\mathfrak{m}},Y\in\mathfrak{m}_2,$ or $[Y,W]_{\mathfrak{m}},X\in\mathfrak{m}_2.$ This occurs precisely when the triple $(X,Y,W)$ is, up to permutation, one of the following:
\[
(Z,A_j,A_k),\ (Z,A_j,B_k),\ (Z,B_j,B_k),\ j,k=1,...,p.
\]
A straightforward computation gives
\[
(d\omega_2)_{\gamma(t)}\left(X^*(\gamma(t)),Y^*(\gamma(t)),W^*(\gamma(t))\right)=0
\]
for all these triples. Consequently, $d\omega_2=0.$\\

The proof that $d\omega_3=0$ is analogous and we omit it.   
\end{proof}

\begin{proposition}\label{prop:kahler-balanced}
The following assertions hold:
\begin{enumerate}
\item[$a)$] $(g,J)$ is K\"ahler if and only if
\begin{equation}\label{eq:Kahler-condition}
\begin{cases}
&\displaystyle\frac{f_2'}{f_2}+\frac{\sin(2t)}{2}\frac{f_0^2}{f_2^2}=0,\\[1.2em]
&\displaystyle\frac{f_3'}{f_3}-\frac{\sin(2t)}{2}\frac{f_0^2}{f_3^2}=0,
\end{cases}\,\,\,\,\,0<t<\frac{\pi}{2}.
\end{equation}

\item[$b)$] $(g,J)$ is balanced if and only if
\begin{equation}\label{eq:balanced-condition}
p\left(\frac{f_2'}{f_2}+\frac{\sin(2t)}{2}\frac{f_0^2}{f_2^2}\right)+q\left(\frac{f_3'}{f_3}-\frac{\sin(2t)}{2}\frac{f_0^2}{f_3^2}\right)=0,\,\,\,\, 0<t<\frac{\pi}{2}.
\end{equation}
\end{enumerate}
\end{proposition}
\begin{proof}
By Proposition~\ref{prop:fundamental-form} and Lemma~\ref{lem:differential-forms},
\begin{align*}
d\omega=&d\left(-\sin(2t)f_0^2\,dt\wedge\zeta+f_2^2\omega_2+f_3^2\omega_3\right)\\
=&(-2\cos(2t)f_0^2-2\sin(2t)f_0f_0')dt\wedge dt\wedge \zeta-\sin(2t)f_0^2d(dt)\wedge\zeta+\sin(2t)f_0^2dt\wedge d\zeta\\
&+2f_2f_2'dt\wedge\omega_2+2f_3f_3'dt\wedge\omega_3\\
=&\sin(2t)f_0^2\,dt\wedge d\zeta
+2f_2f_2'\,dt\wedge\omega_2
+2f_3f_3'\,dt\wedge\omega_3\\
=&\sin(2t)f_0^2\,dt\wedge(\omega_2-\omega_3)+2f_2f_2'\,dt\wedge\omega_2
+2f_3f_3'\,dt\wedge\omega_3\\
=&\left(2f_2f_2'+\sin(2t)f_0^2\right)dt\wedge\omega_2
+\left(2f_3f_3'-\sin(2t)f_0^2\right)dt\wedge\omega_3.
\end{align*}
Since the forms $dt\wedge\omega_2$ and $dt\wedge\omega_3$ are linearly independent, $d\omega=0$ on $(\mathbb{C}P^{p+q+1})_0$ if and only if
\[
2f_2f_2'+\sin(2t)f_0^2=0,\ 2f_3f_3'-\sin(2t)f_0^2=0,\ 0<t<\frac{\pi}{2}.
\]
Dividing the first equation by $2f_2^2$ and the second equation by $2f_3^2,$ we obtain \eqref{eq:Kahler-condition}. This proves $a).$\\

To prove $b)$, recall that 
\[
d\omega^{p+q}=(p+q)d\omega\wedge\omega^{p+q-1}.
\]
Since each term of $d\omega$ contains $dt$, only the terms of $\omega^{p+q-1}$ which do not contain $dt$ contribute to $d\omega\wedge\omega^{p+q-1}.$ Thus,
\begin{align*}
d\omega\wedge\omega^{p+q-1}&=d\omega\wedge\left(f_2^2\omega_2+f_3^2\omega_3\right)^{p+q-1}.
\end{align*}
Note that $\omega_2^{p+1}=0,\ \omega_3^{q+1}=0,$ and $dt\wedge\omega_2^p\wedge\omega_3^q\neq0.$ Therefore,
\begin{align*}
\left(f_2^2\omega_2+f_3^2\omega_3\right)^{p+q-1}&=\sum\limits_{r=0}^{p+q-1}\left(\begin{array}{c}
     p+q-1 \\
     r 
\end{array}\right)f_2^{2r}f_3^{2(p+q-1-r)}\omega_2^r\wedge\omega_3^{p+q-1-r}\\[0.5em]
&=\left(\begin{array}{c}
     p+q-1 \\
     p-1 
\end{array}\right)f_2^{2p-2}f_3^{2q}\omega_2^{p-1}\wedge\omega_3^{q}+\left(\begin{array}{c}
     p+q-1 \\
     q-1 
\end{array}\right)f_2^{2p}f_3^{2q-2}\omega_2^{p}\wedge\omega_3^{q-1}
\end{align*}
Consequently,
\begin{align*}
d\omega\wedge\omega^{p+q-1}=&\left\{\left(\begin{array}{c}
     p+q-1 \\
     p-1 
\end{array}\right)\left(2f_2f_2'+\sin(2t)f_0^2\right)f_2^{2p-2}f_3^{2q}\right.\\[0.5em]
&\left.+\left(\begin{array}{c}
     p+q-1 \\
     q-1 
\end{array}\right)\left(2f_3f_3'-\sin(2t)f_0^2\right)f_2^{2p}f_3^{2q-2}\right\}
dt\wedge\omega_2^p\wedge\omega_3^q.
\end{align*}
Hence, $d\omega^{p+q}=0$ on $(\mathbb{C}P^{p+q+1})_0$ if and only if
\[
\left(\begin{array}{c}
     p+q-1 \\
     p-1 
\end{array}\right)\left(2f_2f_2'+\sin(2t)f_0^2\right)f_2^{2p-2}f_3^{2q}+\left(\begin{array}{c}
     p+q-1 \\
     q-1 
\end{array}\right)\left(2f_3f_3'-\sin(2t)f_0^2\right)f_2^{2p}f_3^{2q-2}=0
\]
for all $t\in(0,\pi/2).$ Multiplying this equation by 
\[
\frac{p!q!}{2(p+q-1)!f_2^{2p}f_3^{2q}},
\]
we obtain \eqref{eq:balanced-condition}. The proof is complete.
\end{proof}
In order to study the equation~\eqref{eq:balanced-condition}, we introduce the variable $x=\sin^2(t),\ 0<t<\pi/2,$ and set $h_j(x):=f_j(t)^2,\ j=0,2,3,\ x\in[0,1].$ The next proposition translates the conditions on the functions $f_j$ into conditions on the functions $h_j,\ j=0,2,3.$
\begin{proposition}\label{prop:smoothness-hj}
    The conditions \eqref{eq:f0-smoothness}, \eqref{eq:smoothness-at-Q0} and \eqref{eq:smoothness-at-Q1} are equivalent to
    \begin{equation}\label{eq:smoothness-hj}
    \begin{cases}
        &h_j\in C^{\infty}([0,1]),\ j=0,2,3,\\
        &h_0(x)>0,\ x\in [0,1],\\
        &h_2(x)>0,\ x\in [0,1),\\
        &h_3(x)>0,\ x\in (0,1],\\
        &h_3(0)=0,\ h_3'(0)=h_0(0),\\
        &h_2(1)=0,\ h_2'(1)=-h_0(1).
    \end{cases}
    \end{equation}
\end{proposition}
\begin{proof} 
   Assume that $f_0,f_2,f_3$ are smooth on $[0,\pi/2]$ and satisfy
\eqref{eq:f0-smoothness}, \eqref{eq:smoothness-at-Q0} and
\eqref{eq:smoothness-at-Q1}. Let
\[
\psi(s):=\sum_{k=1}^{\infty}(-1)^{k+1}\frac{2^{2k-1}}{(2k)!}s^k.
\]
Then $\psi$ is smooth near $0$ and for $s\geq0$, $\psi(s)=\sin^2(\sqrt{s}).$ Moreover, $\psi(0)=0,$ $\psi'(0)=1.$ Hence, $\psi$ has a local smooth inverse $\psi^{-1}$ defined near $0$ and $\psi^{-1}(0)=0,$ $(\psi^{-1})'(0)=1.$ Near $t=0$, by \eqref{eq:f0-smoothness} and \eqref{eq:smoothness-at-Q0}, we have $f_0(t)=\phi_0(t^2),$ $f_2(t)=\phi_{2,0}(t^2),$ $f_3(t)=t\phi_{3,0}(t^2),$ where $\phi_0,\phi_{2,0},\phi_{3,0}$ are smooth near $0$ and $\phi_{3,0}(0)=f_0(0).$ Since $x=\psi(t^2),$ we have $t^2=\psi^{-1}(x)$. Thus $h_0(x)=\phi_0(\psi^{-1}(x))^2,$ $h_2(x)=\phi_{2,0}(\psi^{-1}(x))^2$ and $h_3(x)=\psi^{-1}(x)\phi_{3,0}(\psi^{-1}(x))^2.$
Thus, $h_0,h_2,h_3$ are smooth near $x=0$, $h_3(0)=0$ and $h_3'(0)=(\psi^{-1})'(0)\phi_{3,0}(0)^2=f_0(0)^2=h_0(0).$ Now, let $r=\pi/2-t.$ Near $r=0$, by \eqref{eq:f0-smoothness} and \eqref{eq:smoothness-at-Q1}, $f_0\left(\pi/2-r\right)=\phi_1(r^2),$ $f_2\left(\pi/2-r\right)=r\phi_{2,1}(r^2),$ $f_3\left(\pi/2-r\right)=\phi_{3,1}(r^2),$ where $\phi_{2,1}(0)=f_0\left(\pi/2\right).$ Since $1-x=\cos^2t=\sin^2r=\psi(r^2),$ we have $r^2=\psi^{-1}(1-x)$. Hence, $h_0(x)=\phi_1(\psi^{-1}(1-x))^2,$ $h_2(x)=\psi^{-1}(1-x)\phi_{2,1}(\psi^{-1}(1-x))^2$ and $h_3(x)=\phi_{3,1}(\psi^{-1}(1-x))^2.$ Thus, $h_0,h_2,h_3$ are smooth near $x=1,$ $h_2(1)=0$ and $h_2'(1)=-(\psi^{-1})'(0)\phi_{2,1}(0)^2=-f_0\left(\pi/2\right)^2=-h_0(1).$ Finally, since the functions $f_j,\ j=0,2,3$ are positive and smooth on $(0,\pi/2),$ the functions $h_j,\ j=0,2,3$ are positive and smooth on $(0,1)$ and $h_0(0)=f_0(0)^2>0,$ $h_0(1)=f_0(\pi/2)^2>0,$ $h_2(0)=f_2(0)^2>0,$ $h_3(1)=f_3(\pi/2)^2>0.$ Hence, the functions $h_j,\ j=0,2,3$ satisfy \eqref{eq:smoothness-hj}.\\

Conversely, assume that the functions $h_j$, $j=0,2,3$, satisfy \eqref{eq:smoothness-hj}. Since $h_j(\sin^2(t))=f_j(t)^2,\ j=0,2,3,$ we have $f_j(t)=\sqrt{h_j(\sin^2(t))},\ j=0,2,3.$ Also, $h_0$ is smooth and positive on $[0,1]$, so the functions $\phi_0(s):=\sqrt{h_0(\psi(s))},$ $\phi_1(s):=\sqrt{h_0(1-\psi(s))}$ are smooth and positive near $0,$ and we have
\[
f_0(t)=\phi_0(t^2),\ f_0\left(\frac{\pi}{2}-t\right)=\phi_1(t^2).
\]
Similarly, since $h_2$ is smooth and positive near $x=0$, the function $\phi_{2,0}(s):=\sqrt{h_2(\psi(s))}$ is smooth and positive near $0$ and we have $f_2(t)=\phi_{2,0}(t^2).$ From $h_2(1)=0,$ we obtain that there exists a smooth function $\widehat{h}_2$ near $1$ such that $h_2(x)=(1-x)\widehat{h}_2(x).$ Furthermore, $\widehat{h}_2(1)=-h_2'(1)=h_0(1)>0,$ so $\widehat h_2$ is positive near $1$. Thus, 
\[
f_2\left(\pi/2-t\right)=\sqrt{h_2(\cos^2(t))}=\sin(t)\,\sqrt{\widehat{h}_2(\cos^2(t))}=t\phi_{2,1}(t^2),
\]
where 
\[
\phi_{2,1}(s):=\frac{\sin(\sqrt{s})}{\sqrt{s}}\sqrt{\widehat{h}_2(1-\psi(s))}
\]
is smooth and positive near $0$ and $\phi_{2,1}(0)=\sqrt{\widehat{h}_2(1)}=\sqrt{h_0(1)}=f_0(\pi/2).$\\

Analogously, since $h_3$ is smooth and positive near 1, the function $\phi_{3,1}(s):=\sqrt{h_3(1-\psi(s))}$ is smooth and positive near $0$ and we have $f_3(\pi/2-t)=\phi_{3,1}(t^2).$ On the other hand, since $h_3(0)=0$, there exists a smooth function $\widehat{h}_3$ near $0$ such that $h_3(x)=x\widehat{h}_3(x).$ This function satisfies $\widehat{h}_3(0)=h_3'(0)=h_0(0)>0,$ so $\widehat{h}_3$ is positive near $0$. It follows that
\[
f_3(t)=\sqrt{h_3(\sin^2t)}=\sin(t)\,\sqrt{\widehat{h}_3(\sin^2(t))}=t\phi_{3,0}(t^2),
\]
where 
\[\phi_{3,0}(s)=\frac{\sin(\sqrt{s})}{\sqrt{s}}\sqrt{\widehat{h}_3(\psi(s))}
\]
is smooth and positive near $0$ and $\phi_{3,0}(0)=\sqrt{\widehat{h}_3(0)}=\sqrt{h_0(0)}=f_0(0).$ Hence,
\eqref{eq:f0-smoothness}, \eqref{eq:smoothness-at-Q0} and
\eqref{eq:smoothness-at-Q1} hold and the proof is complete.
\end{proof}
Note that the equations \eqref{eq:Kahler-condition} and \eqref{eq:balanced-condition} are equivalent to
\begin{equation}\label{eq:Kahler-condition-hj}
\begin{cases}
&\displaystyle h_2'+h_0=0,\\
&\displaystyle h_3'-h_0=0,
\end{cases}\,\,\,\,\,0<x<1,
\end{equation}
and
\begin{equation}\label{eq:balanced-condition-hj}
   p\frac{h_2'+h_0}{h_2}+q\frac{h_3'-h_0}{h_3}=0,\ 0<x<1, 
\end{equation}
respectively. We will study the latter equation under the conditions~\eqref{eq:smoothness-hj}.
\begin{proposition}\label{prop:solutions-balanced}
Let $(h_0,h_2,h_3)$ be a solution of \eqref{eq:balanced-condition-hj} satisfying \eqref{eq:smoothness-hj}. Then there exists a unique function $V\in C^\infty([0,1])$ such that $V(0)=0$ and
\begin{equation}\label{eq:h2-h3-formulas}
\begin{cases}
&\displaystyle h_2(x)=e^{qV(x)}\int_x^1 h_0(s)e^{-qV(s)}\,ds\\[1.2em]
&\displaystyle h_3(x)=e^{-pV(x)}\int_0^x h_0(s)e^{pV(s)}\,ds.
\end{cases}
\end{equation}

Conversely, let $h_0\in C^\infty([0,1])$ be positive, let $V\in C^\infty([0,1])$ with $V(0)=0$ and define $h_2,h_3$ by the formulas \eqref{eq:h2-h3-formulas}. Then $(h_0,h_2,h_3)$ is a solution of \eqref{eq:balanced-condition-hj} satisfying \eqref{eq:smoothness-hj}. Moreover, $(h_0,h_2,h_3)$ satisfies \eqref{eq:Kahler-condition-hj} if and only if $V=0.$
\end{proposition}
\begin{proof} Assume that $(h_0,h_2,h_3)$ is a solution of \eqref{eq:balanced-condition-hj} that satisfies \eqref{eq:smoothness-hj}. Then,
\[
\frac{1}{q}\frac{h_2'+h_0}{h_2}=-\frac{1}{p}\frac{h_3'-h_0}{h_3}=:W(x),\ 0<x<1.
\]
Since $h_j\in C^{\infty}([0,1]),\ j=0,2,3,$ $h_2>0$ on $[0,1)$ and $h_3>0$ on $(0,1],$ $W$ admits a unique smooth extension to $[0,1].$ Define
\[
V(x):=\int_{0}^xW(y)\,dy,\ 0\leq x\leq 1.
\]
Then $V\in C^{\infty}([0,1]),$ $V(0)=0$ and
\[
\begin{cases}
    &h_2'-qV'h_2=-h_0\\
    &h_3'+pV'h_3=h_0.
\end{cases}
\]
Since $h_2(1)=0$ and $h_3(0)=0,$ solving these two linear equations yields
\[
h_2(x)=e^{qV(x)}\int_x^1 h_0(s)e^{-qV(s)}\,ds
\]
and
\[
h_3(x)=e^{-pV(x)}\int_0^x h_0(s)e^{pV(s)}\,ds,
\]
for all $x\in[0,1].$ If $\widetilde{V}\in C^{\infty}([0,1])$ is such that $\widetilde{V}(0)=0$ and the formulas above hold with $\widetilde{V}$ instead of $V,$ then $\widetilde{V}'(x)=W(x)=V'(x),\ 0\leq x\leq 1.$ Since $\widetilde{V}(0)=0=V(0),$ $\widetilde{V}=V.$ Thus, $V$ is unique.\\

Conversely, let $h_0\in C^\infty([0,1])$ be positive, let $V\in C^\infty([0,1])$ with $V(0)=0$ and define $h_2, h_3$ by the formulas \eqref{eq:h2-h3-formulas}. Clearly, $h_2,h_3\in C^\infty([0,1])$. Since $h_0$ is positive, we have $h_2(x)>0,$ $0\leq x<1,$ $h_3(x)>0,\ 0<x\leq1,$ $h_2(1)=0,\ h_3(0)=0.$ Differentiating the expressions for $h_2$ and $h_3,$ we obtain
\[
h_2'=qV'h_2-h_0,\ h_3'=-pV'h_3+h_0.
\]
Hence, $h_2'(1)=-h_0(1),$ $h_3'(0)=h_0(0)$ and $(h_0,h_2,h_3)$ satisfies \eqref{eq:smoothness-hj}. Furthermore,
\[
\frac{h_2'+h_0}{h_2}=qV',\ 
\frac{h_3'-h_0}{h_3}=-pV',
\]
on $(0,1)$. Thus
\[
p\frac{h_2'+h_0}{h_2}+q\frac{h_3'-h_0}{h_3}=pqV'-pqV'=0,\ 0<x<1,
\]
and $(h_0,h_2,h_3)$ is a solution of \eqref{eq:balanced-condition-hj}.\\

Finally, by the identities
\[
\frac{h_2'+h_0}{h_2}=qV',\ \frac{h_3'-h_0}{h_3}=-pV',
\]
the system \eqref{eq:Kahler-condition-hj} holds if and only if $V'=0$ on $(0,1)$. Since $V(0)=0$, this is equivalent to $V=0$ on $[0,1].$
\end{proof}
As an immediate consequence of Corollary~\ref{coro:compatible-metrics}, Proposition~\ref{prop:smooth-extension}, Proposition~\ref{prop:kahler-balanced}, Proposition~\ref{prop:smoothness-hj} and Proposition~\ref{prop:solutions-balanced}, we have the following result. 
\begin{corollary}\label{coro:parametrization-balanced-metrics} The following assertions hold:
\begin{itemize}
    \item[$a)$] The set of $G_{p,q}$-invariant balanced Hermitian structures $(g,J)$ on $\mathbb{C}P^{p+q+1}$ is in bijection with the set
    \[
    \mathcal{SK}:=\left\{(h_0,V)\in C^{\infty}([0,1])\times C^{\infty}([0,1]):h_0>0,\ V(0)=0\right\}.
    \]
    \item[$b)$] The set of $G_{p,q}$-invariant K\"ahler structures  $(g,J)$ on $\mathbb{C}P^{p+q+1}$ is in bijection with the set
    \[
    \mathcal{K}:=\left\{h_0\in C^{\infty}([0,1]):h_0>0\right\}.
    \]
\end{itemize}
\end{corollary}
\section{An application: inhomogeneous metrics with positive scalar curvature}\label{section:application}
In this section, we construct an explicit one-parameter family of $G_{m,m}$-invariant balanced Hermitian metrics on $\mathbb{C}P^{2m+1}$ and study the scalar curvature of this family.
\begin{proposition}\label{prop:explicit-family}
Let $m\geq2.$ For $\varepsilon\in\mathbb{R}$, consider $h_{0,\varepsilon}(x)=1,$ $V_\varepsilon(x)=\varepsilon x,$
and let
\[
h_{2,\varepsilon}(x)=
\begin{cases}
\displaystyle\frac{1-e^{-m\varepsilon(1-x)}}{m\varepsilon},& \textnormal{if}\ \varepsilon\neq0,\\[1.2em]
1-x,& \textnormal{if}\ \varepsilon=0,
\end{cases}\,\,\,\,\,\,\, h_{3,\varepsilon}(x)=
\begin{cases}
\displaystyle\frac{1-e^{-m\varepsilon x}}{m\varepsilon},& \textnormal{if}\ \varepsilon\neq0,\\[1.2em]
x,& \textnormal{if}\ \varepsilon=0.
\end{cases}
\]
Let $g_\varepsilon=(1,h_{2,\varepsilon},h_{3,\varepsilon})$ be the corresponding $G_{m,m}$-invariant Hermitian metric on $\mathbb{C}P^{2m+1}$. Then:
\begin{itemize}
    \item[$a)$] $(g_\varepsilon,J)$ is balanced for every $\varepsilon\in\mathbb{R}.$
    \item[$b)$] $g_0=g_{\mathrm{FS}}$ and $(g_\varepsilon,J)$ is non-K\"ahler for every $\varepsilon\neq0.$
    \item[$c)$] There exists $\varepsilon_0>0$ such that, for every $0<|\varepsilon|<\varepsilon_0$, the metric $g_\varepsilon$ is inhomogeneous and has positive scalar curvature.
\end{itemize}
\end{proposition}
\begin{proof}
For $p=q=m$, the formulas in Proposition~\ref{prop:solutions-balanced} become
\[
h_2(x)=e^{mV(x)}\int_x^1 h_0(s)e^{-mV(s)}\,ds,\ h_3(x)=e^{-mV(x)}\int_0^x h_0(s)e^{mV(s)}\,ds.
\]
Taking $h_0=1$ and $V=V_\varepsilon=\varepsilon x$, for $\varepsilon\neq0$ we obtain
\[
h_{2,\varepsilon}(x)=e^{m\varepsilon x}\int_x^1e^{-m\varepsilon s}\,ds=\frac{1-e^{-m\varepsilon(1-x)}}{m\varepsilon}
\]
and
\[
h_{3,\varepsilon}(x)=e^{-m\varepsilon x}\int_0^xe^{m\varepsilon s}\,ds=\frac{1-e^{-m\varepsilon x}}{m\varepsilon}.
\]
For $\varepsilon=0$, we have $h_{2,0}(x)=1-x,$ $h_{3,0}(x)=x.$ It follows from Corollary~\ref{coro:parametrization-balanced-metrics} that $(g_\varepsilon,J)$ is balanced for every $\varepsilon\in\mathbb{R}.$ This proves $a).$\\

Now, for all $t\in(0,\pi/2),$  we have
\begin{align*}
    &\sqrt{h_0(\sin^2(t))}=1\\
    &\sqrt{h_{2,0}(\sin^2(t))}=\sqrt{1-\sin^2(t)}=\cos(t)\\
    &\sqrt{h_{3,0}(\sin^2(t))}=\sqrt{\sin^2(t)}=\sin(t).
\end{align*}
Thus $g_0=(1,h_{2,0},h_{3,0})$ corresponds to the metric $(1,\sin(2t),\cos(t),\sin(t),0)$ which, by Example~\ref{ex:Fubini-Study}, is the Fubini-Study metric. Furthermore, if $\varepsilon\neq0$, then $V_\varepsilon\neq0.$ Hence, by Proposition~\ref{prop:solutions-balanced}, $(g_\varepsilon,J)$ is non-K\"ahler. This proves $b).$\\

To prove $c)$, observe that, since $\mathbb{C}P^{2m+1}$ is compact and $g_{\mathrm{FS}}$ has positive scalar curvature, it would be sufficient to show that $g_\varepsilon\to g_0=g_{\mathrm{FS}}$ in the $C^\infty$-topology to obtain the positivity of the scalar curvature of $g_\varepsilon$ for sufficiently small $\varepsilon$. However, we will follow a different approach that will also allow us to prove that $g_\varepsilon$ is inhomogeneous for sufficiently small nonzero $\varepsilon$. We will provide an explicit formula for the scalar curvature $s_\varepsilon$ of $g_\varepsilon$. Since $g_\varepsilon$ is $G_{m,m}$-invariant, it suffices to compute $s_\varepsilon$ at points $\gamma(t),$ $t\in[0,\pi/2].$ Let $\varepsilon\neq 0$ and let
\begin{align*}
    f_0(t)=\sqrt{h_0(\sin^2(t))}=1,\ f_{2,\varepsilon}(t)=\sqrt{h_{2,\varepsilon}(\sin^2(t))},\ f_{3,\varepsilon}(t)=\sqrt{h_{3,\varepsilon}(\sin^2(t))}.
\end{align*}
Thus $g_\varepsilon=(1,\sin(2t),f_{2,\varepsilon},f_{3,\varepsilon},0).$ By Lemma~\ref{lem:arc-length-Q0}, the curve $\gamma\big{|}_{(0,\pi/2)}$ is a unit-speed geodesic with respect to $g_{\varepsilon},$ orthogonal to all principal orbits. Therefore, we can apply the formulas in \cite[Proposition~1.14]{GZ2002} for the Ricci curvature of $g_\varepsilon$ to obtain 
\begin{align*}
&\mathrm{Ric}_{\gamma(t)}(\gamma'(t),\gamma'(t))=4
-2m\frac{f_{2,\varepsilon}''(t)}{f_{2,\varepsilon}(t)}
-2m\frac{f_{3,\varepsilon}''(t)}{f_{3,\varepsilon}(t)},\\[1.2em]
&\begin{aligned}
\mathrm{Ric}_{\gamma(t)}\left(\frac{Z^*(\gamma(t))}{\sin(2t)},\frac{Z^*(\gamma(t))}{\sin(2t)}\right)=&\,4+\frac{m}{2}\frac{\sin^2(2t)}{f_{2,\varepsilon}(t)^4}+\frac{m}{2}\frac{\sin^2(2t)}{f_{3,\varepsilon}(t)^4}\\[1.2em]
&-4m\frac{\cos(2t)f_{2,\varepsilon}'(t)}{\sin(2t)f_{2,\varepsilon}(t)}-4m\frac{\cos(2t)f_{3,\varepsilon}'(t)}{\sin(2t)f_{3,\varepsilon}(t)},
\end{aligned}\\[1.2em]
&\begin{aligned}
\mathrm{Ric}_{\gamma(t)}\left(\frac{X^*(\gamma(t))}{f_{2,\varepsilon}(t)},\frac{X^*(\gamma(t))}{f_{2,\varepsilon}(t)}\right)=&-\frac{f_{2,\varepsilon}''(t)}{f_{2,\varepsilon}(t)}+\frac{2(m+1)}{f_{2,\varepsilon}(t)^2}-\frac{1}{2}\frac{\sin^2(2t)}{f_{2,\varepsilon}(t)^4}\\[1.2em]
&-(2m-1)\left(\frac{f_{2,\varepsilon}'(t)}{f_{2,\varepsilon}(t)}\right)^2-2\frac{\cos(2t)f_{2,\varepsilon}'(t)}{\sin(2t)f_{2,\varepsilon}(t)}-2m\frac{f_{2,\varepsilon}'(t)f_{3,\varepsilon}'(t)}{f_{2,\varepsilon}(t)f_{3,\varepsilon}(t)},
\end{aligned}\\[1.2em]
&\begin{aligned}
\mathrm{Ric}_{\gamma(t)}\left(\frac{Y^*(\gamma(t))}{f_{3,\varepsilon}(t)},\frac{Y^*(\gamma(t))}{f_{3,\varepsilon}(t)}\right)=&-\frac{f_{3,\varepsilon}''(t)}{f_{3,\varepsilon}(t)}+\frac{2(m+1)}{f_{3,\varepsilon}(t)^2}-\frac{1}{2}\frac{\sin^2(2t)}{f_{3,\varepsilon}(t)^4}\\[1.2em]
&-(2m-1)\left(\frac{f_{3,\varepsilon}'(t)}{f_{3,\varepsilon}(t)}\right)^2-2\frac{\cos(2t)f_{3,\varepsilon}'(t)}{\sin(2t)f_{3,\varepsilon}(t)}-2m\frac{f_{2,\varepsilon}'(t)f_{3,\varepsilon}'(t)}{f_{2,\varepsilon}(t)f_{3,\varepsilon}(t)},
\end{aligned}
\end{align*}
where $0<t<\pi/2,$ $X\in\{A_j,B_j:j=1,...,m\},$ and $Y\in\{C_j,D_j:j=1,...,m\}.$ Hence,
\begin{align*}
    s_\varepsilon(t)=&\,\mathrm{Ric}_{\gamma(t)}(\gamma'(t),\gamma'(t))+\mathrm{Ric}_{\gamma(t)}\left(\frac{Z^*(\gamma(t))}{\sin(2t)},\frac{Z^*(\gamma(t))}{\sin(2t)}\right)\\[1.2em]
    &+\sum\limits_{j=1}^m\left(\mathrm{Ric}_{\gamma(t)}\left(\frac{A_j^*(\gamma(t))}{f_{2,\varepsilon}(t)},\frac{A_j^*(\gamma(t))}{f_{2,\varepsilon}(t)}\right)+\mathrm{Ric}_{\gamma(t)}\left(\frac{B_j^*(\gamma(t))}{f_{2,\varepsilon}(t)},\frac{B_j^*(\gamma(t))}{f_{2,\varepsilon}(t)}\right)\right)\\[1.2em]
    &+\sum\limits_{j=1}^m\left(\mathrm{Ric}_{\gamma(t)}\left(\frac{C_j^*(\gamma(t))}{f_{3,\varepsilon}(t)},\frac{C_j^*(\gamma(t))}{f_{3,\varepsilon}(t)}\right)+\mathrm{Ric}_{\gamma(t)}\left(\frac{D_j^*(\gamma(t))}{f_{3,\varepsilon}(t)},\frac{D_j^*(\gamma(t))}{f_{3,\varepsilon}(t)}\right)\right)\\[1.2em]
   =&\,8+4m(m+1)\left(\frac{1}{f_{2,\varepsilon}(t)^2}+\frac{1}{f_{3,\varepsilon}(t)^2}\right)-\frac{m}{2}\sin^2(2t)\left(\frac{1}{f_{2,\varepsilon}(t)^4}+\frac{1}{f_{3,\varepsilon}(t)^4}\right)\\[1.2em]
   &-8m\frac{\cos(2t)}{\sin(2t)}\left(\frac{f_{2,\varepsilon}'(t)}{f_{2,\varepsilon}(t)}+\frac{f_{3,\varepsilon}'(t)}{f_{3,\varepsilon}(t)}\right)-4m\left(\frac{f_{2,\varepsilon}''(t)}{f_{2,\varepsilon}(t)}+\frac{f_{3,\varepsilon}''(t)}{f_{3,\varepsilon}(t)}\right)\\[1.2em]
   &-2m(2m-1)\left(\left(\frac{f_{2,\varepsilon}'(t)}{f_{2,\varepsilon}(t)}\right)^2+\left(\frac{f_{3,\varepsilon}'(t)}{f_{3,\varepsilon}(t)}\right)^2\right)-8m^2\frac{f_{2,\varepsilon}'(t)f_{3,\varepsilon}'(t)}{f_{2,\varepsilon}(t)f_{3,\varepsilon}(t)}.
\end{align*}
In terms of $x=\sin^2(t),$ we have
\begin{align*}
    &\sin(2t)=2\sqrt{x(1-x)},\ \cos(2t)=1-2x,\\[0.5em]
    &f_{2,\varepsilon}(t)=\sqrt{h_{2,\varepsilon}(x)},\ f_{3,\varepsilon}(t)=\sqrt{h_{3,\varepsilon}(x)},\\[0.5em]
    &\frac{f_{2,\varepsilon}'(t)}{f_{2,\varepsilon}(t)}=\sqrt{x(1-x)}\frac{h_{2,\varepsilon}'(x)}{h_{2,\varepsilon}(x)},\ \frac{f_{3,\varepsilon}'(t)}{f_{3,\varepsilon}(t)}=\sqrt{x(1-x)}\frac{h_{3,\varepsilon}'(x)}{h_{3,\varepsilon}(x)},\\[0.5em]
    &\frac{f_{2,\varepsilon}''(t)}{f_{2,\varepsilon}(t)}=(1-2x)\frac{h_{2,\varepsilon}'(x)}{h_{2,\varepsilon}(x)}-x(1-x)\left(\frac{h_{2,\varepsilon}'(x)}{h_{2,\varepsilon}(x)}\right)^2+2x(1-x)\frac{h_{2,\varepsilon}''(x)}{h_{2,\varepsilon}(x)},\\[0.5em]
    &\frac{f_{3,\varepsilon}''(t)}{f_{3,\varepsilon}(t)}=(1-2x)\frac{h_{3,\varepsilon}'(x)}{h_{3,\varepsilon}(x)}-x(1-x)\left(\frac{h_{3,\varepsilon}'(x)}{h_{3,\varepsilon}(x)}\right)^2+2x(1-x)\frac{h_{3,\varepsilon}''(x)}{h_{3,\varepsilon}(x)},
\end{align*}
therefore
\begin{align*}
    s_\varepsilon(x)=&\,8+4m(m+1)\left(\frac{1}{h_{2,\varepsilon}(x)}+\frac{1}{h_{3,\varepsilon}(x)}\right)-2mx(1-x)\left(\frac{1}{h_{2,\varepsilon}(x)^2}+\frac{1}{h_{3,\varepsilon}(x)^2}\right)\\[1.2em]
    &-8m(1-2x)\left(\frac{h_{2,\varepsilon}'(x)}{h_{2,\varepsilon}(x)}+\frac{h_{3,\varepsilon}'(x)}{h_{3,\varepsilon}(x)}\right)-8mx(1-x)\left(\frac{h_{2,\varepsilon}''(x)}{h_{2,\varepsilon}(x)}+\frac{h_{3,\varepsilon}''(x)}{h_{3,\varepsilon}(x)}\right)\\[1.2em]
    &-2m(2m-3)x(1-x)\left(\left(\frac{h_{2,\varepsilon}'(x)}{h_{2,\varepsilon}(x)}\right)^2+\left(\frac{h_{3,\varepsilon}'(x)}{h_{3,\varepsilon}(x)}\right)^2\right)-8m^2x(1-x)\frac{h_{2,\varepsilon}'(x)h_{3,\varepsilon}'(x)}{h_{2,\varepsilon}(x)h_{3,\varepsilon}(x)}.
\end{align*}
Since we are assuming $\varepsilon\neq 0,$ 
\begin{align*}
&h_{2,\varepsilon}'(x)=m\varepsilon h_{2,\varepsilon}(x)-1,\ h_{2,\varepsilon}''(x)=m^2\varepsilon^2h_{2,\varepsilon}(x)-m\varepsilon,\\
&h_{3,\varepsilon}'(x)=1-m\varepsilon h_{3,\varepsilon}(x),\ h_{3,\varepsilon}''(x)=-m\varepsilon +m^2\varepsilon^2h_{3,\varepsilon}(x).
\end{align*}
Thus,
\begin{align*}
    s_\varepsilon(x)=&\,8-4m^3\varepsilon^2x(1-x)+4m(m+1)\left(\frac{1}{h_{2,\varepsilon}(x)}+\frac{1}{h_{3,\varepsilon}(x)}\right)+8m(1-2x)\left(\frac{1}{h_{2,\varepsilon}(x)}-\frac{1}{h_{3,\varepsilon}(x)}\right)\\[1.2em]
    &-4m(m-1)x(1-x)\left(\frac{1}{h_{2,\varepsilon}(x)^2}+\frac{1}{h_{3,\varepsilon}(x)^2}\right)-4m^2\varepsilon x(1-x)\left(\frac{1}{h_{2,\varepsilon}(x)}+\frac{1}{h_{3,\varepsilon}(x)}\right)\\[1.2em]
    &+8m^2x(1-x)\left(\frac{1}{h_{2,\varepsilon}(x)}\frac{1}{h_{3,\varepsilon}(x)}\right).
\end{align*}
Consider the function
\[
G(u):=\begin{cases}
    \displaystyle\frac{1-e^{-u}}{u},&\textnormal{if}\ u\neq 0,\\[1.2em]
    1, &\textnormal{if}\ u=0.
\end{cases}
\]
This function is analytic because
\[
G(u)=\sum\limits_{k=0}^\infty\frac{(-1)^{k}}{(k+1)!}u^{k},\ u\in\mathbb{R},
\]
and it satisfies $G(0)=1,\ G'(0)=-1/2,\ G''(0)=1/3.$ Therefore, the function $F(u):=1/G(u)$ is also analytic in a neighborhood of zero and 
\[
F(0)=1,\ F'(0)=-\frac{G'(0)}{G(0)^2}=\frac{1}{2},\ F''(0)=2\frac{G'(0)^2}{G(0)^3}-\frac{G''(0)}{G(0)^2}=\frac{1}{6},
\]
so we can write
\[
F(u)=1+\frac{1}{2}u+\frac{1}{12}u^2+R(u),
\]
where there exist $C>0$ and $\delta>0$ such that $|R(u)|\leq C|u|^3$ whenever $|u|<\delta.$\\

For all $x\in[0,1)$,
\begin{align*}
\frac{1}{h_{2,\varepsilon}(x)}&=\frac{m\varepsilon}{1-e^{-m\varepsilon(1-x)}}\\
&=\frac{1}{1-x}F\left(m\varepsilon(1-x)\right)\\
&=\frac{1}{1-x}+\frac{m}{2}\varepsilon+\frac{m^2}{12}(1-x)\varepsilon^2+\frac{R\left(m\varepsilon(1-x)\right)}{1-x}\\
&=\frac{1}{1-x}+\frac{m}{2}\varepsilon+\frac{m^2}{12}(1-x)\varepsilon^2+R_2(x,\varepsilon),
\end{align*}
where
\[
R_2(x,\varepsilon):=\frac{R(m\varepsilon(1-x))}{1-x}.
\]
Note that, if $0<|\varepsilon|<\delta/m,$ then $|m\varepsilon(1-x)|<\delta,$ for all $x\in[0,1).$ Consequently,
\[
|R_2(x,\varepsilon)|=\left|\frac{R\left(m\varepsilon(1-x)\right)}{1-x}\right|\leq Cm^3|\varepsilon|^3(1-x)^2\leq Cm^3|\varepsilon|^3,
\]
and we can write
\[
\frac{1}{h_{2,\varepsilon}(x)}=\frac{1}{1-x}+\frac{m}{2}\varepsilon+\frac{m^2}{12}(1-x)\varepsilon^2+O(\varepsilon^3(1-x)^2). 
\]
Analogously, for every $x\in(0,1]$,
\begin{align*}
\frac{1}{h_{3,\varepsilon}(x)}&=\frac{m\varepsilon}{1-e^{-m\varepsilon x}}\\
&=\frac{1}{x}F\left(m\varepsilon x\right)\\
&=\frac{1}{x}+\frac{m}{2}\varepsilon+\frac{m^2}{12}x\varepsilon^2+\frac{R\left(m\varepsilon x\right)}{x}\\
&=\frac{1}{x}+\frac{m}{2}\varepsilon+\frac{m^2}{12}x\varepsilon^2+R_3(x,\varepsilon),
\end{align*}
where
\[
R_{3}(x,\varepsilon):=\frac{R(m\varepsilon x)}{x}.
\]
As before, $0<|\varepsilon|<\delta/m$ implies $|R_3(x,\varepsilon)|\leq Cm^3|\varepsilon|^3x^2\leq Cm^3|\varepsilon|^3.$ Thus,
\[
\frac{1}{h_{3,\varepsilon}(x)}=\frac{1}{x}+\frac{m}{2}\varepsilon+\frac{m^2}{12}x\varepsilon^2+O(\varepsilon^3x^2). 
\]
For $x\in(0,1),$ it follows that
\begin{align*}
&\frac{1}{h_{2,\varepsilon}(x)}+\frac{1}{h_{3,\varepsilon}(x)}=\frac{1}{x(1-x)}+m\varepsilon
+\frac{m^2}{12}\varepsilon^2+O(\varepsilon^3),\\[0.5em]
&\frac{1}{h_{2,\varepsilon}(x)}-\frac{1}{h_{3,\varepsilon}(x)}=\frac{2x-1}{x(1-x)}+\frac{m^2}{12}(1-2x)\varepsilon^2+O(\varepsilon^3),\\[0.5em]
&\begin{aligned}\frac{1}{h_{2,\varepsilon}(x)^2}=&\frac{1}{(1-x)^2}+\frac{m}{1-x}\varepsilon+\frac{5m^2}{12}\varepsilon^2+\left(\frac{1}{1-x}O(\varepsilon^3(1-x)^2)+\frac{m^3}{12}(1-x)\varepsilon^3\right),\\[0.5em]
=&\frac{1}{(1-x)^2}+\frac{m}{1-x}\varepsilon+\frac{5m^2}{12}\varepsilon^2+O(\varepsilon^3)
\end{aligned}\\[0.5em]
&\begin{aligned}
\frac{1}{h_{3,\varepsilon}(x)^2}=&\frac{1}{x^2}+\frac{m}{x}\varepsilon+\frac{5m^2}{12}\varepsilon^2+\left(\frac{1}{x}O(\varepsilon^3x^2)+\frac{m^3}{12}x\varepsilon^3\right)\\[0.5em]
=&\frac{1}{x^2}+\frac{m}{x}\varepsilon+\frac{5m^2}{12}\varepsilon^2+O(\varepsilon^3),
\end{aligned}
\end{align*}
where the term $O(\varepsilon^3)$ is uniform in $x\in(0,1).$ Consequently,
\[
x(1-x)\left(\frac{1}{h_{2,\varepsilon}(x)^2}+\frac{1}{h_{3,\varepsilon}(x)^2}\right)=\frac{x^2+(1-x)^2}{x(1-x)}+m\varepsilon+\frac{5m^2}{6}x(1-x)\varepsilon^2
+O(\varepsilon^3)
\]
and
\begin{align*}
\frac{x(1-x)}{h_{2,\varepsilon}(x)h_{3,\varepsilon}(x)}=&1+\frac{m}{2}\varepsilon+\frac{m^2}{12}\bigl(1+x(1-x)\bigr)\varepsilon^2+\left(xO(\varepsilon^3)+(1-x)O(\varepsilon^3)+\frac{m^3}{24}x(1-x)\varepsilon^3\right)\\[0.5em]
=&1+\frac{m}{2}\varepsilon+\frac{m^2}{12}\bigl(1+x(1-x)\bigr)\varepsilon^2+O(\varepsilon^3).
\end{align*}
Substituting these expansions into the expression for $s_\varepsilon(x)$, we obtain
\begin{align*}
s_\varepsilon(x)=&\,8-4m^3\varepsilon^2x(1-x)+4m(m+1)\left(\frac{1}{x(1-x)}+m\varepsilon+\frac{m^2}{12}\varepsilon^2+O(\varepsilon^3)\right)\\[0.5em]
&+8m(1-2x)\left(\frac{2x-1}{x(1-x)}+\frac{m^2}{12}(1-2x)\varepsilon^2
+O(\varepsilon^3)\right)\\[0.5em]
&-4m(m-1)\left(\frac{x^2+(1-x)^2}{x(1-x)}+m\varepsilon+\frac{5m^2}{6}x(1-x)\varepsilon^2+O(\varepsilon^3)\right)\\[0.5em]
&-4m^2\varepsilon\left(1+mx(1-x)\varepsilon+\frac{m^2}{12}x(1-x)\varepsilon^2+O(\varepsilon^3)\right)\\[0.5em]
&+8m^2\left(1+\frac{m}{2}\varepsilon+\frac{m^2}{12}\bigl(1+x(1-x)\bigr)\varepsilon^2+O(\varepsilon^3)\right)\\[0.5em]
=&A_0(x)+A_1(x)\varepsilon+A_2(x)\varepsilon^2+O(\varepsilon^3),
\end{align*}
where
\begin{align*}
    &\begin{aligned}
        A_0(x)=&8+8m^2+\frac{4m(m+1)}{x(1-x)}-\frac{8m(1-2x)^2}{x(1-x)}-\frac{4m(m-1)(x^2+(1-x)^2)}{x(1-x)}\\[0.5em]
        =&8(2m+1)(m+1),
    \end{aligned}\\[0.5em]
    &\begin{aligned}
        A_1(x)=&4m^2(m+1)-4m^2(m-1)-4m^2+4m^3\\[0.5em]
        =&4m^2(m+1),
    \end{aligned}\\[0.5em]
    &\begin{aligned}
    A_2(x)=&-4m^3x(1-x)+\frac{m^3(m+1)}{3}+\frac{2m^3}{3}(1-2x)^2\\[0.5em]
    &-\frac{10m^3(m-1)}{3}x(1-x)-4m^3x(1-x)+\frac{2m^4}{3}(1+x(1-x))\\
    =&m^3\left(m+1-\frac{8m+22}{3}x(1-x)\right).
\end{aligned}
\end{align*}
Therefore,
\[
s_\varepsilon(x)=8(m+1)(2m+1)+4m^2(m+1)\varepsilon+m^3\left(m+1-\frac{8m+22}{3}x(1-x)\right)\varepsilon^2+O(\varepsilon^3).
\]
Since $x(1-x)\in[0,1/4]$ for all $x\in(0,1),$ we obtain
\[
s_\varepsilon(x)=8(m+1)(2m+1)+O(\varepsilon)
\]
uniformly in $x\in(0,1).$ In particular,
\[
\lim\limits_{\varepsilon\to 0}\sup\limits_{x\in(0,1)}|s_\varepsilon(x)-8(m+1)(2m+1)|=0.
\]
Hence, there exists $\varepsilon_0>0$ such that $0<|\varepsilon|<\varepsilon_0$ implies
\[
|s_\varepsilon(x)-8(m+1)(2m+1)|<4(m+1)(2m+1),\ \forall x\in(0,1).
\]
Consequently, if $0<|\varepsilon|<\varepsilon_0,$ then
\[
s_\varepsilon(x)>4(m+1)(2m+1)>0,\ \forall x\in(0,1).
\]
By continuity of the scalar curvature we have $s_\varepsilon(0),s_\varepsilon(1)\geq4(m+1)(2m+1)>0.$ That is, $s_\varepsilon>0$ for $0<|\varepsilon|<\varepsilon_0.$ Finally, to prove that $g_\varepsilon$ is inhomogeneous, note that
\begin{align*}
s_\varepsilon\left(\frac14\right)-s_\varepsilon\left(\frac12\right)=\frac{m^3(4m+11)}{24}\varepsilon^2+O(\varepsilon^3).
\end{align*}
Since $m^3(4m+11)/24>0,$ we may decrease $\varepsilon_0,$ if necessary, so that
\[
s_\varepsilon\left(\frac14\right)\neq s_\varepsilon\left(\frac12\right)
\]
for $0<|\varepsilon|<\varepsilon_0$. Thus, $g_\varepsilon$ has non-constant scalar curvature and is therefore inhomogeneous. The proof is complete.
\end{proof}

\section*{Acknowledgements}
L. G. is partially supported by The São Paulo Research Foundation (FAPESP) grants 2023/13131-8 and 2024/00923-6; and by National Council for Scientific and Technological Development (CNPq) grant 306021/2024-2.

\section*{Ethical Statement}

\paragraph{Data Availability.}
No datasets were generated or analyzed during the current study. Therefore, data sharing is not applicable to this article.
\\

\paragraph{Conflict of Interest.}
The authors declare that they have no conflict of interest.
\\

\paragraph{Research Involving Human Participants and/or Animals.}
This article does not contain any studies involving human participants or animals performed by any of the authors.

\bibliographystyle{acm}

\begin{thebibliography}{99}

\bibitem{AS2021}
Alonso, I., and Salvatore, F. On the existence of balanced metrics on six-manifolds of cohomogeneity one. {\it Ann. Glob. Anal. Geom.} 61(2) (2022): 309-331.

\bibitem{AS2024}
Alonso, I., and Salvatore, F. Correction to: On the existence of balanced metrics on six-manifolds of cohomogeneity one. {\it Ann. Glob. Anal. Geom.} 67(1) (2025): Article 3.

\bibitem{AV2016}
Andrada, A., and Villacampa, R. Abelian balanced Hermitian structures on unimodular Lie algebras. {\it Transform. Groups} 21(4) (2016): 903-927.

\bibitem{AP2022}
Angella, D., and Pediconi, F. On cohomogeneity one Hermitian non-K\"ahler metrics. {\it Proc. R. Soc. Edinb., Sect. A} 153(2) (2023): 545-587.

\bibitem{AGG2026}
Araujo, A. L. A., Grajales, B., and Grama, L. Cohomogeneity one Einstein metrics on complex projective spaces. arXiv:2605.28539 (2026).

\bibitem{B2026}
Broder, K. Cohomogeneity-one Einstein metrics on closed $4$-manifolds. arXiv:2605.24859 (2026).

\bibitem{CR2026}
Cao-Labora, G., and Rodr\'iguez-V\'azquez, A. Inhomogeneous Einstein metrics on complex projective spaces. arXiv:2608.16880 (2026).

\bibitem{C2023}
Correa, E. M. $t$-Gauduchon Ricci-flat metrics on non-K\"ahler Calabi-Yau manifolds. arXiv:2310.01788 (2023).


\bibitem{F2024}
Fino, A. Canonical metrics in complex geometry. {\it Boll. Unione Mat. Ital.} 18(1) (2025): 185-198.

\bibitem{FGV2019}
Fino, A., Grantcharov, G., and Vezzoni, L. Astheno-K\"ahler and balanced structures on fibrations. {\it Int. Math. Res. Not. IMRN} 2019(22) (2019): 7093-7117.

\bibitem{FP2023}
Fino, A., and Paradiso, F. Balanced Hermitian structures on almost abelian Lie algebras. {\it J. Pure Appl. Algebra} 227(2) (2023): Article 107186.

\bibitem{FLY2012}
Fu, J., Li, J., and Yau, S.-T. Balanced metrics on non-K\"ahler Calabi-Yau threefolds. {\it J. Differ. Geom.} 90(1) (2012): 81-129.

\bibitem{GZ2002}
Grove, K., and Ziller, W. Cohomogeneity one manifolds with positive Ricci curvature. {\it Invent. Math.} 149(3) (2002): 619-646.


\bibitem{L2026}
Lauret, J. Compact homogeneous complex manifolds. arXiv:2608.30747 (2026).

\bibitem{M1982}
Michelsohn, M. L. On the existence of special metrics in complex geometry. {\it Acta Math.} 149(1) (1982): 261-295.

\bibitem{T1973}
Takagi, R. On homogeneous real hypersurfaces in a complex projective space. {\it Osaka J. Math.} 10(3) (1973): 495-506.


\bibitem{VZ}
Verdiani, L., and Ziller, W. Smoothness conditions in cohomogeneity one manifolds. {\it Transform. Groups} 27(1) (2022): 311-342.

\bibitem{Whitney1943}
Whitney, H. Differentiable even functions. {\it Duke Math. J.} 10(1) (1943): 159-160.
\end{thebibliography}

\end{document}